\documentclass[10pt]{amsart}
\usepackage{amsfonts}
\usepackage{color}
\usepackage[square,compress,comma, numbers]{natbib}
\usepackage[colorlinks=true, citecolor=blue, linkcolor=blue]{hyperref}
\allowdisplaybreaks[4]
\usepackage{amssymb}
\usepackage{amsmath}
\usepackage{ulem}

\definecolor{c20}{rgb}{0.,0.7,0.}
\definecolor{c30}{rgb}{0.,0.,1.}
\definecolor{c40}{rgb}{1,0.1,0.7}
\definecolor{c50}{rgb}{1,0,0}
\definecolor{c60}{rgb}{1,0.9,0.1}

\allowdisplaybreaks[4]

\def\cLa#1{\textcolor{black}{#1}}

\def\bl#1{\textcolor{black}{#1}}

\newcommand{\kb}[1]{\boldsymbol{#1}}
\newcommand{\vk}[1]{\kb{#1}}

\newcommand{\abs}[1]{\left\lvert #1 \right\rvert}

\newcommand{\EF}[1]{\mathbb{E}\left\{#1\right\}}

\newcommand{\pk}[1]{\mathbb{P} \left \{#1 \right \} }

\newcommand{\R}{\mathbb{R}}

\newcommand{\N}{\mathbb{N}}
\newcommand{\inr}{\in \R}

\newcommand{\ldot}{,\ldots,}

\newcommand{\BQN}{\begin{eqnarray}}
\newcommand{\EQN}{\end{eqnarray}}
\newcommand{\BQNY}{\begin{eqnarray*}}
\newcommand{\EQNY}{\end{eqnarray*}}

\newcommand{\BS}{\begin{sat}}
\newcommand{\ES}{\end{sat}}
\newcommand{\BT}{\begin{theo}}
\newcommand{\ET}{\end{theo}}
\newcommand{\BL}{\begin{lem}}
\newcommand{\EL}{\end{lem}}
\newcommand{\BK}{\begin{korr}}
\newcommand{\EK}{\end{korr}}

\newcommand{\BD}{\begin{de}}
\newcommand{\ED}{\end{de}}

\newcommand{\BIT}{\begin{itemize}}
\newcommand{\EIT}{\end{itemize}}
\newcommand{\BDI}{\begin{description}}
\newcommand{\EDI}{\end{description}}

\newcommand{\BRM}{\begin{remarks}}
\newcommand{\ERM}{\end{remarks}}

\newcommand{\BEL}{\begin{lem}}
\newcommand{\EEL}{\end{lem}}

\newtheorem{theo}{Theorem}[section]
\newtheorem{sat}[theo]{Proposition}
\newtheorem{de}[theo]{Definition}
\newtheorem{lem}[theo]{Lemma}

\newtheorem{example}[theo]{Example}
\newtheorem{korr}[theo]{Corollary}
\newtheorem{remark}[theo]{Remark}
\newtheorem{remarks}[theo]{Remarks}
\newtheorem{prop}[theo]{Proposition}

\newcommand{\nelem}[1]{{Lemma \ref{#1}}}

\newcommand{\netheo}[1]{{Theorem \ref{#1}}}
\newcommand{\nekorr}[1]{{Corollary \ref{#1}}}

\newcommand{\prooftheo}[1]{ \textsc{\bf Proof of Theorem} \ref{#1}:}
\newcommand{\proofprop}[1]{\textsc{\bf Proof of Proposition} \ref{#1}:}
\newcommand{\prooflem}[1]{\textsc{\bf Proof of Lemma} \ref{#1}:}
\newcommand{\proofkorr}[1]{\textsc{\bf Proof of Corollary} \ref{#1}:}

\newcommand{\COM}[1]{}

\newcommand{\QED}{\hfill $\Box$}

\def\rw{\rightarrow}

\def\IF{\infty}

\def\oo{(1+o(1))}

\def\LT{\left}
\def\RT{\right}

\def\rw{\rightarrow}

\def\vn{\varepsilon}
\def\Var{\text{Var}}
\def\E{E(u)}

\def\II{\mathbb{I}}

\def\Bu+#1{\mathcal{B}^{\varepsilon+}_{u}(#1)}

\begin{document}

\title[Extremes of  vector-valued Gaussian processes with Trend]{Extremes of  vector-valued Gaussian processes with Trend}

\author{Long Bai}
\address{Long Bai,
Department of Actuarial Science, %\\Faculty of Business and Economics\\
University of Lausanne\\
UNIL-Dorigny, 1015 Lausanne, Switzerland
}
\email{Long.Bai@unil.ch}

\author{Krzysztof D\c{e}bicki}
\address{Krzysztof D\c{e}bicki, Mathematical Institute, University of Wroc\l aw, pl. Grunwaldzki 2/4, 50-384 Wroc\l aw, Poland}
\email{Krzysztof.Debicki@math.uni.wroc.pl}

%\author{Enkelejd  Hashorva}
%\address{Enkelejd Hashorva, Department of Actuarial Science, %\\Faculty of Business and Economics\\
%University of Lausanne,\\
%UNIL-Dorigny, 1015 Lausanne, Switzerland
%}
%\email{Enkelejd.Hashorva@unil.ch}

\author{Peng Liu}
\address{Liu, Peng, Department of Actuarial Science, University of Lausanne, UNIL-Dorigny, 1015 Lausanne, Switzerland
}
\email{peng.liu@unil.ch}

\bigskip

\date{\today}
 \maketitle

{\bf Abstract:} Let $\vk{X}(t)=(X_1(t), \dots, X_n(t)), t\in \mathcal{T}\subset \R $ be a centered vector-valued Gaussian process with
independent components and continuous trajectories, and $\vk{h}(t)=(h_1(t),\dots, h_n(t)), t\in \mathcal{T} $
be a vector-valued continuous function. We investigate the asymptotics of
$$\pk{\sup_{t\in \mathcal{T} } \min_{1\leq i\leq n}(X_i(t)+h_i(t))>u}$$
as $u\rw\IF$.
As an illustration to the derived results we analyze two important classes of $\vk{X}(t)$:  with
locally-stationary structure and with varying variances of the coordinates,
and calculate exact asymptotics of simultaneous ruin
probability and ruin time in a Gaussian risk model.
%Additionally we find exact asymptotics
%of ruin probability in Gaussian
%two : the variance of $X_i$ is non-constant and $X_i$ is locally stationary. Our results have potential application in risk theory to describe simultaneous ruin probability for multiple dependent companies.



\section{Introduction and Preliminaries}
Motivated by various applied-oriented problems,
the asymptotics of
\BQN\label{onedimensional}
\pk{\sup_{t\in \mathcal{T}}(X(t)+h(t))>u},
\EQN
as $u\rw\IF$,
for both $\mathcal{T}=[0,T]$ and $\mathcal{T}=[0,\IF)$,
where $X(t)$ is a centered Gaussian process with continuous trajectories and
$h(t)$ is a continuous function, attracted substantial interest in the literature; see e.g.
\cite{HP99, DE2002, HP2004, DI2005,HA2013, KEP2016, KP2015, KP2017}
and references therein
%For example,
%ruin probability
%one of the classical problems in risk theory, related with ruin probability
%
%is $X$ has stationary increments and $h(t)=-ct$, $c>0$,  (\ref{onedimensional})
%gives the probability or ruin in a Gaussian risk model,
%
%One of important areas
for
connections of (\ref{onedimensional}) with problems considered, e.g., in risk theory or  fluid queueing models.
For example, in the setting
of risk theory
 one usually supposes that  $h(t)=-ct$, with $c>0$ and
$X$ has stationary increments.
Then, using that
$\pk{\sup_{t\in \mathcal{T}}(X(t)+h(t))>u}=\pk{\inf_{t\in \mathcal{T}}(u-X(t)+ct)<0}$,
(\ref{onedimensional}) represents {\it ruin probability}, with
$X(t)$ modelling
the accumulated claims amount in time interval $[0,t]$,
$c$ being the constant premium rate and $u$, the initial capital.
The most celebrated model in this context is the Brownian risk model
introduced in
the seminal work by Iglehart \cite{Igl69}, where $X$ is a standard Brownian motion.
Extensions to more general class of Gaussian processes
with stationary increments,
including fractional Brownian motions, was analyzed in, e.g., \cite{Mic98,HP99,HP2004,Zhang081,Zhang082}.
Recent interest in the analysis of risk models has turned to
the investigation of multidimensional ruin problems, including investigation of {\it simultaneous ruin }
probability of some number, say $n$, of independent risk processes
\[
\pk{\exists_{t\in \mathcal{T}}\forall_{i=1,...,n} (u_i-X_i(t)+c_it)<0},
\]
see, e.g., \cite{Pal082} and \cite{Pal081}. Motivated by this sort of problems,
in this paper we investigate multidimensional counterpart of (\ref{onedimensional}), i.e.,
we are interested in the exact asymptotics of
\BQN\label{multidimensional}
\pk{\exists_{t\in[0,T]}\vk{X}(t)+\vk{h}(t)>u\vk{1}}=
\pk{\sup_{t\in [0,T] } \min_{1\leq i\leq n}(X_i(t)+h_i(t))>u},
\EQN
as $u\rw\IF$, $T\in(0,\infty)$, where
$\vk{X}(t)=(X_1(t), \dots, X_n(t)), t\in \mathcal{T}\subset \R $ is an $n-$dimensional centered Gaussian process with mutually
independent coordinates and continuous trajectories
and
$\vk{h}(t)=(h_1(t),\dots, h_n(t)), t\in [0,T] $ is a vector-valued continuous function.

We note that (\ref{multidimensional}) can also be viewed as the probability that the conjunction set
$\mathcal{S}_{T,u}:=\{t\in[0,T]:\min _{1\le i\le n} (X_i(t)+h_i(t))>u\}$
is not empty in Gaussian
conjunction problem, since
\[
\pk{\mathcal{S}_{T,u}\neq \emptyset}=\pk{\sup_{t\in [0,T] } \min_{1\leq i\leq n}(X_i(t)+h_i(t))>u},
\]
see, e.g., \cite{Wro00,Tabis} and references therein.

The main results of this contribution extend recent findings of \cite{Tabis}, where
the exact asymptotics of (\ref{multidimensional}) for $h_i\equiv0, 1\leq i\leq n$ %and $\mathcal{T}=[0,T]$
was analyzed; see also  \cite{DHJ16}
where $\vk{X}(t)$ is a multidimensional Brownian motion, $h_i(t)=c_it$ and $T=\infty$,
and \cite{DKMR10,PS05} for LDP-type results.
 It appears that the presence of the drift function substantially increases difficulty of the problem
when comparing it with the analysis given for the driftless case in \cite{Tabis}.
More specifically, as advocated in Section 2, it
requires to
deal with
\BQN
\pk{\sup_{t\in [0,T] } \min_{1\leq i\leq n}X_{u,i}(t)>u},\nonumber
\EQN
where
$(X_{u,i}(t), t\in [0,T])_u$, $i=1,...,n$ are families (with respect to $u$) of centered threshold-dependent Gaussian processes;
see Theorem \ref{PreThm1}.
%which leads to qualitatively new types of the asymptotics.

In Section \ref{s.3} we apply general results derived in Section 2 to two important
families of Gaussian processes, i.e. i) to locally-stationary processes in the sense of Berman and
ii) to processes with varying variance $\Var(X_i(t))$, $t\in [0,T]$.
Then, as an example to the derived theory, we
analyze the probability of simultaneous
ruin in Gaussian risk model.
Complementary,
we investigate the limit distribution of the {\it simultaneous ruin time}
$$\tau_u:=\inf\{t\geq 0: \LT(\vk{X}(t)+\vk{h}(t)\RT) >u\vk{1}\},$$
conditioned that $\tau_u\le T$, as $u\to\infty$.

%A relate to  \cite{DHJ16}  for a work that deals with a special instance of (\ref{multidimensional}),
%where $\vk{X}(t)$ is a multidimensional Brownian motion and $h_i(t)=c_it$.
%Logarithmic asymptotics of (\ref{multidimensional}) under general assumptions
%on dependence structure of $\vk{X}(t)$  and $\mathcal{T}\subset \mathbb{R}^d$ was derived in, see also \cite{PS05}.

%In comparison to \cite{Tabis},
%where the centered vector-valued Gaussian processes are investigated,
%it turns out that the trend function $\vk{h}$ makes the problem much more difficult in the sense that we need to consider the extremes of the threshold-dependent
%centered Gaussian vector processes $\vk{X}_u(t), u>0$.
%Thus, in section 3 we analyze extremes of $\vk{X}_u(t), u>0$ over short interval. T
%hen applying this main results to the cases that the variance of $X_i$ is non-constant and $X_i$  is locally stationary,
%we present our results of vector-valued processes with trend in Section 4.\\

Organization of the rest of the paper:  Section 2 is devoted to the main result of this contribution, concerning the extremes of
the threshold-dependent centered Gaussian vector processes. In Section 3 we specify our result to
locally-stationary vector-valued Gaussian processes with trend and
non-stationary Gaussian vector-valued processes with trend. Detailed proofs of all the results are postponed to Section 4.
Additionally, in Section 3 we analyze asymptotics of the simultaneous ruin probability.

\COM{
Let $\{\vk{X}(t),t\geq 0\}$ be a vector-valued Gaussian process where $\vk{X}(t)=(X_1(t),\ldots,X_n(t))$ with $X_i(t),t\geq 0, i=1,\ldots,n, n\in\N$, being independent centered Gaussian processes with almost surely (a.s.) continuous sample paths. The asymptotic behavior of the probability that $\vk{X}(t)$ enters the upper orthant $\{(x_1,\ldots,x_n):x_i>u,i=1,\ldots,n\}$ over a fixed time interval $[0,T]$, i.e.,
\BQN\label{p1}
p_{T,u}:=\pk{\exists_{t\in[0,T]}\forall_{i=1,\ldots,n}X_i(t)>u}=\pk{\sup_{t\in[0,T]}
\min_{1\leq i\leq n}X_i(t)>u},\ u\rw\IF,
\EQN
was firstly investigated in \cite{MR1747100} in the case of smooth Gaussian random fields, further contributions can be found in \cite{MR2775212,ChengXiao13}. Results for non-Gaussian random fields and general stationary processes can be found in \cite{MR2654766,DebOrderStats}.\\
In \cite{Tabis, Anshin05,Debicki10}, the exact asymptotics of $p_{T,u}$ with non-stationary Gaussian processes  $X_i's$ was derived. But some drawbacks are exists in the literature, such as  $X_i's$ are centered in \cite{Tabis}. We further generalize these processes in this paper.\\
Another motivation is from \cite{GauTrend16}, and there asymptotics of a general class of  non-centered Gaussian process with trend, i.e.,
\BQN
\pk{\sup_{t\in[0,T]}X(t)+h(t)>u},\ u\rw\IF,
\EQN
 is considered where $X(t)$ is a non-stationary centered Gaussian process and $h(t)$ is a bounded function. In this paper this will be a special case of the Gaussian processes conjunction problem \eqref{p1}.\\

}

%\section{Supremum over Threshold-dependent Intervals} %Preliminaries
%\section{Extremes of vector-valued Gaussian processes with trend}
\section{Main Results}\label{s.2}
We begin with observation that,
for sufficiently large $u$,
\begin{eqnarray}
\pk{\sup_{t\in [0,T] } \min_{1\leq i\leq n}(X_i(t)+h_i(t))>u}=
\pk{\exists_{t\in[0,T]}\vk{X}_u (t)>u\vk{1}},
\end{eqnarray}
where
$\vk{X}_u (t)=\left(\frac{uX_1(t)}{u-h_1(t)}, \dots, \frac{uX_n(t)}{u-h_n(t)}\right)$
is a family of centered vector-valued threshold-dependent Gaussian processes.
Since the above rearrangement appears to be  useful for the technique of the proof
that we use in order to get the exact asymptotics of (\ref{multidimensional}), then
in this section we focus on asymptotics of extremes of threshold-dependent vector-valued
Gaussian processes.

%In the section we  analyse
%the tail probability of a general centered vector-valued Gaussian process
%over a threshold-dependent interval, which plays
%the key role in the investigation of the exact asymptotics of vector-valued Gaussian processes with trend. Let ${\bf X}(t):=(X_1(t), \dots, X_n(t))$ with $X_i, 1\leq i\leq n$ being independent Gaussian processes with mean zero and continuous trajectories and $\vk{h}(t)$ being a vector function.  We are interest in the asymptotics, as $u\rw\IF$, of
%$$\pk{\exists_{t\in[0,T]}\LT(\vk{X}(t)+{\vk{h}}(t)\RT) >u\vk{1}},$$
%with $\vk{c}$ being a constant vector.
%Rewrite the above probability as
%\BQNY
%\pk{\exists_{t\in[0,T]}\LT(\vk{X}(t)+{\vk{h}}(t)\RT) >u\vk{1}}=\pk{\exists_{t\in[0,T]}\left(\frac{uX_1(t)}{u-h_1(t)}, \dots, \frac{uX_n(t)}{u-h_n(t)}\right) >u\vk{1}}.
%\EQNY
%This indicates that in order to investigate the extremes of vector-valued Gaussian processes with trend, it suffices to deal with the extremes of threshold-dependent centered vector-valued Gaussian processes.

More specifically, let ${\bf X}_u(t):=(X_{u,1}(t),\dots, X_{u,n}(t)), t\in E(u)$,
with $0\in E(u)=(x_1(u), x_2(u)),$
be a family of centered $n$-dimensional vector-valued Gaussian processes
with continuous trajectories.
Let $\sigma^2_{u,i}(\cdot)$ and $r_{u,i}(\cdot,\cdot)$
be the variance function
and the correlation function
of $X_{u,i}(t)$, $1\leq i\leq n$ respectively.
Moreover, we tacitly
assume that $X_{u,i}(t)$, $1\leq i\leq n$ are
mutually independent.

We shall impose the following assumptions on ${\bf X}_u(t)$:\\
% For any large $u$, the following generalized variance function
%\BQNY
%g_u(t)=\sum_{i=1}^{n}\frac{1}{\sigma_{u,i}^2(t)}
%\EQNY
%attains its minimum over $E_u$  at the unique point $0$ and}
\textbf{A1}:
$
\lim_{u\rw\IF}\vk{\sigma}_u(0)=\vk{\sigma}>\vk{0}.
$

\textbf{A2}: There exist $\lambda_{i}\in[0,\IF), 1\leq i\leq n$ with $\max_{1\leq i\leq n}\lambda_{i}>0$
and some continuous functions
$f_{i}(\cdot),1\leq i\leq n$ with  $f_i(0)=0$
such that
for any $\epsilon\in (0,1)$, as $u\rightarrow\IF$,
$$
\left|\left(\frac{\sigma_{u,i}(0)}{\sigma_{u,i}(t)}-1\right)u^2-f_{i}
 (u^{\lambda_{i}}t)\right|\leq  \epsilon(\LT|f_{i}(u^{\lambda_{i}}t)\RT|+1), \quad t\in\E.
$$
\textbf{A3}: There exist $\alpha_i\in(0,2]$ and $\ a_i>0,\ 1\leq i\leq n$ such that
$$
\lim_{u\rightarrow\IF}\underset{t\not=s }{\sup_{s,t\in\E}}\left|\frac{1-r_{u,i}(t,s)}{a_i|t-s|^{\alpha_i}}-1\right|=0.
$$

\COM{In \netheo{PreThm1}, we shall investigate the asymptotic of
\BQN
\pk{\exists_{t\in\E} \vk{X}_{u}(t+t_u)>\vk{M}_u},\ u\rw\IF.
\EQN
for $\vk{M}_u$ satisfying $\lim_{u\rw\IF}\frac{\vk{M}_u}{u}=c\vk{1}$ with $c>0$.}
\COM{For a  closed interval $E \subset \R$, denote by  $C^*_0(E)$ the collection of all  functions satisfying
that  $f(\cdot)$ is continuous over $E$, $f(0)=0$, and there exist $\epsilon_2>\epsilon_1>0$ such that
$$\lim_{|t|\rw\IF,t\in E}|f(t)|/|t|^{\epsilon_1}=\IF, \lim_{|t|\rw\IF,t\in E}|f(t)|/|t|^{\epsilon_2}=0,$$
provided that  $\sup\{x:x\in E\}=\IF$ or $\inf\{x:x\in E\}=-\IF$.\\}

In the following we write $f\in\mathcal{R}_{\alpha}$ to denote that function $f$ is
regularly varying at $\IF$ with index $\alpha$, see \cite{EKM97,Res,Soulier}
for the definition and properties of regularly varying functions.\\
Let  $\lambda:=\max_{1\leq i\leq n}\lambda_{i}$,
$\alpha:=\min_{1\leq i\leq n}\alpha_i$,
$\cLa{\widetilde{\vk{f}}(t):=\LT(\widetilde{f}_1(t),\ldots,\widetilde{f}_n(t)\RT)}$
with \cLa{$$\widetilde{f}_i\LT(t\RT)=f_i\LT(t\RT)\mathbb{I}_{\{\lambda_i=\lambda\}}$$} and suppose that
$x_1(u)\in \mathcal{R}_{-\mu_1},\ x_2(u)\in \mathcal{R}_{-\mu_2}$ with $\mu_1,\mu_2\geq \lambda$ and
\BQN\label{xx}
&&\lim_{u\rw\IF}u^{\lambda}x_1(u)=x_1\in[-\IF,\IF),\nonumber\\
&&\lim_{u\rw\IF}u^{\lambda}x_2(u)=x_2\in(-\IF,\IF], \ x_1<x_2,\\
&&\lim_{u\rw\IF}u^{\lambda_j}x_i(u)=0, i=1,2, \lambda_j<\lambda.\nonumber
\EQN
If $|x_1|+|x_2|=\IF$, we additionally
assume that
\BQN\label{FF}
%\sum_{i=1}^n\frac{\widehat{f}_i(t)}{\sigma_i^2}\not\equiv 0, t\in[x_1,x_2],
%\ \ \ \
\underset{t\in [x_1,x_2]}{\liminf_{|t|\rw\IF}}\LT(\sum_{i=1}^n
\frac{\widetilde{f}_i(t)}{\sigma_i^2}\RT)\Big/
{\LT(\sum_{i=1}^n\frac{\abs{\widetilde{f}_i(t)}}{\sigma_i^2}\RT)}>0.
\EQN
Assumption
(\ref{FF}) means that the negative components of $\frac{\widetilde{f}_i(t)}{\sigma_i^2}, 1\leq i\leq n$
do not play a significant role to the sum in comparison with the positive components.

Moreover, we suppose that $0\cdot\IF =0$, $u^{-\IF}=0$ for any $u>0$ and introduce
$$[x_1,x_2]:=\lim_{u\rw\IF}f(u)[x_1(u),x_2(u)],$$
if $\lim_{u\rw\IF}f(u)x_1(u)=x_1\in[-\IF,\IF)$ and $\lim_{u\rw\IF}f(u)x_2(u)=x_2\in(-\IF,\IF]$ with $x_1<x_2$.

Next we introduce some notation and definition of the Pickands-Piterbarg constants.

% that appear
%in all the obtained in next sections asymptotics.
Throughout this paper, all the operations on vectors are meant componentwise, for instance, for any given
 $ \vk{x} = (x_1,\ldots,x_n)\in \R ^n$ and $\vk{y} = (y_1,\ldots,y_n) \in \R ^n $, we write $ \vk{x} > \vk{y} $ if and only if
  $ x_i > y_i $  for all $ 1 \leq i \leq n $, {write $1/\vk{x}=(1/x_1,\cdots,1/x_n)$ if $x_i\neq 0, 1 \leq i \leq n$}, and write $\vk{x}\vk{y}=(x_1y_1 \ldot x_ny_n)$.
Further we set  $ \vk{0}  := (0,\ldots,0)\in\R^n $ and $ \vk{1} : = (1,\ldots,1)\in\R^n$. \\
\COM{We shall refer to $\{\vk{X}_u(t), t\ge0\}$
as a family of centered $n$-dimensional  vector-valued Gaussian processes,
where $ \vk{X}_u(t)=(X_{u,1}(t) \ldot X_{u,n}(t))$ with $X_{u,i}$'s being independent centered Gaussian processes with a.s. continuous sample paths.
Since $n$ hereafter is always fixed we shall occasionally omit
"$n$-dimensional", mentioning simply that $\vk{X}_u$ is a centered vector-valued Gaussian process.}
\COM{Let $f$ be a continuous function such that
$$\lim_{|t|\rw\IF,t\in F}|f(t)|/|t|^{\epsilon_1}=\IF, \lim_{|t|\rw\IF,t\in F}|f(t)|/|t|^{\epsilon_2}=0.$$}
%The following are the vector version Pickands and Piterbarg constants.\\
Define for $S_1,S_2 \in \R, \ S_1<S_2$, $\vk{a}=(a_1,a_2,\ldots, a_n)$ with $a_i\geq0,\ 1\leq i\leq n$ and $\vk{f}(t)=(f_1(t),\ldots, f_n(t))$ with $f_i(t),\ 1\leq i\leq n$ being continuous functions
\BQNY
\mathcal{P}_{\alpha,\vk{a}}^{\vk{f}}[S_1,S_2]&:=&\int_{\R^n}e^{\sum_{i=1}^n w_i}\pk{\exists_{t\in[S_1,S_2]}\LT(\sqrt{2\vk{a}}{\vk{B}}_\alpha(t)-{\vk{a}}|t|^\alpha
-{\vk{f}}(t)\RT)>{\vk{w}}}d{\vk{w}}\\
&=&\int_{\R^n}e^{\sum_{i=1}^n w_i}\pk{\sup_{t\in[S_1,S_2]}\LT(\min_{1\leq i\leq n}\sqrt{2a_i}{B}_{\alpha,i}(t)-a_i|t|^\alpha -f_i(t)-w_i\RT)>0}d{\vk{w}}\in(0,\IF),
\EQNY
where  ${\vk{B}}_\alpha(t),t\in\mathbb{R}$ is an $n$-dimensional vector-valued standard fractional Brownian motion (fBm)  with mutually independent coordinates $B_{\alpha,i}(t)$ and common Hurst index $\alpha/2\in(0,1]$. Let
\BQNY
\mathcal{P}_{\alpha,\vk{a}}^{\vk{f}}[0,\IF):=\lim_{S_2\rw\IF} \mathcal{P}_{\alpha,\vk{a}}^{\vk{f}}[0,S_2],\ \
\mathcal{P}_{\alpha,\vk{a}}^{\vk{f}}(-\IF,\IF):=\lim_{S_1\rw-\IF,S_2\rw\IF} \mathcal{P}_{\alpha,\vk{a}}^{\vk{f}}[S_1,S_2].
\EQNY
Let, for $\vk{a}>0$, %Further $\mathcal{H}_{\alpha,\vk{a}}$ defined by
$$\mathcal{H}_{\alpha,\vk{a}}=\lim_{T\rightarrow\IF}\frac{1}{T}\mathcal{P}_{\alpha,\vk{a}}^{\vk{0}}[0,T].
$$
Finiteness of $\mathcal{H}_{\alpha,\vk{a}}$, $\mathcal{P}_{\alpha,\vk{a}}^{\vk{f}}[0,\IF)$
and $\mathcal{P}_{\alpha,\vk{a}}^{\vk{f}}(-\IF,\IF)$
is guaranteed under some restrictions on $\vk{f}(\cdot)$ which are satisfied in
our setup; see \cite{Tabis,GauTrend16,GeneralPit16}.
We refer to, e.g.,
\cite{GeneralPit16,PicandsA,Pit72, DE2002,DI2005,DE2014,DiekerY,DEJ14,Pit20, DM, SBK, Htilt,DHL14Ann,DebKo2013}
for properties of the above constants.

Throughout this paper we write $ f(u) = h(u)(1+o(1))$ or $f(u) \sim h(u)$ if $ \lim_{u \to \infty} \frac{f(u)}{h(u)} = 1 $
and $ f(u) = o(h(u)) $ if $ \lim_{u \to \infty} \frac{f(u)}{h(u)} = 0 $.
Let $\Psi(\cdot)$  denote  the tail distribution of an
$N(0,1)$ random variable,
$\Gamma(\cdot)$ denote the Euler Gamma function
and $\mathbf{I}_{\{\vk{a}=\vk{b}\}}:=(\II_{\{a_1=b_1\}},\ldots,\II_{\{a_n=b_n\}})$ with $\II_{\{\cdot\}}$ being the indicator function.

\BT\label{PreThm1}
Let ${\bf X}_u(t),t\in E(u)$ be a family of centered vector-valued  Gaussian processes
with continuous trajectories and independent coordinates  satisfying {\bf A1-A3} and (\ref{xx})-(\ref{FF}) holds.
Let further  $\cLa{\vk{m}_u}$  be a vector function of $u$ with $\lim_{u\rw\IF}\frac{\vk{m}_u}{u}=\vk{1}$ and for $j\in\{1\leq i\leq n: \lambda_i=\lambda\}$, $f_j(t)$ be regularly varying at $\pm\IF$ with positive index.
Then we have
\BQNY
\pk{\exists_{t\in\E}\vk{X}_{u}(t)>\vk{m}_u}\sim u^{(\frac{2}{\alpha}-\lambda)_{+}}\prod_{i=1}^n
\Psi\LT(\frac{m_{u,i}}{\sigma_{u,i}(0)}\RT)
\times\left\{
\begin{array}{ll}	\mathcal{H}_{\alpha,\frac{\vk{a}}{{\vk{\sigma}}^2}
\mathbf{I}_{\{\vk{\alpha}=\alpha\vk{1}\}}}	\int_{x_1}^{x_2}e^{-\sum_{i=1}^n\frac{\widetilde{f}_i(t)}{\sigma_i^2}}dt, &\hbox{if} \ \ \lambda<2/\alpha,\\	\mathcal{P}^{\frac{\widetilde{\vk{f}}}{\vk{\sigma}^2}}_{\alpha,\frac{\vk{a}}{{\vk{\sigma}}^2}
\mathbf{I}_{\{\vk{\alpha}=\alpha\vk{1}\}}}[x_1,x_2],& \hbox{if} \ \ \lambda=2/\alpha,\\
\int_{\R^n}e^{\sum_{i=1}^n w_i}\mathbb{I}_{\LT\{\exists_{t\in[x_1,x_2]}
-\frac{\widetilde{\vk{f}}(t)}{\vk{\sigma}^2}>\vk{w}\RT\}}d\vk{w},&  \hbox{if} \ \ \lambda>2/\alpha.
\end{array}
\right.
\EQNY
\ET
%\begin{remark} i)\\
%ii)
%\netheo{PreThm1} covers also the case where
%$\lim_{u\rw\IF}\frac{\vk{M}_u}{u}=\vk{c}$ with $c_i>0, 1\leq i\leq n$, since
%$$\pk{\exists_{t\in\E}\vk{X}_{u}(t)>\vk{M}_u}
%=\pk{\exists_{t\in\E}\frac{\vk{X}_{u}(t)}{\vk{c}}>\frac{\vk{M}_u}{\vk{c}}}, \quad \lim_{u\rw\IF}\frac{\vk{M}_u}{u\vk{c}}=\vk{1}.$$
%\end{remark}

\section{Applications}\label{s.3}
In this section we apply \netheo{PreThm1} to the analysis of the exact asymptotics of
\[
\pk{\exists_{t\in[0,T]}\LT(\vk{X}(t)+{\vk{h}}(t)\RT) >u\vk{1}},
\]
as $u\to\infty$.
We distinguish two classes of processes $\vk{X}$: processes with
non-stationary coordinates and processes
with locally-stationary coordinates, including strictly stationary case.

\subsection{Non-stationary coordinates}

Let ${\bf X}(t),t\geq 0$ be a centered vector-valued Gaussian process with
independent coordinates.
% $X_{i}(t)$'s which have continuous trajectories,  variance functions
%$\sigma^2_{i}(\cdot), 1\leq i\leq n$ and correlation functions
%$r_{i}(\cdot,\cdot), 1\leq i\leq n$.  Moreover, assume
Suppose that
$\sigma_i(\cdot), 1\leq i\leq n$
attains its maximum  on $[0, T]$
at the unique point $t_0\in[0,T]$, and further
\BQN\label{eq:sigt0}
\sigma_i(t)=\sigma_i(t_0)-b_i|t-t_0|^{\beta_i}(1+o(1)),\ \
\ \ t\rw t_0
\EQN
with $b_i>0,\beta_i>0$,  and
\BQN\label{eq:rst}
r_i(s,t)=1-a_i|t-s|^{\alpha_i}(1+o(1)),\ \ \ \ s, t \rw t_0
\EQN
for some constants $a_i>0$ and $\alpha_i\in(0,2].$ We further assume that there exists $\mu_1>0$ such that
\BQN\label{holderx}
\max_{i=1,\dots, n}\sup_{s\neq t, s,t\in [0,T]}\frac{\mathbb{E}\left(\left(X_i(t)-X_i(s)\right)^2\right)}{|t-s|^{\mu_1}}<\IF.
\EQN
 Let $\vk{h}(t)$ be a continuous vector function over $[0,T]$ satisfying
 \BQN \label{eq:gtt0}
h_i(t)=h_i(t_0)-c_i|t-t_0|^{\gamma_i}(1+o(1)),\ \ \ \  t \rw t_0
\EQN
with   $c_i<0$ and $\gamma_i\geq\frac{\beta_i}{2}$; and $ c_i\geq0$ and $\gamma_i>0$. Moreover, there exists $\mu_2>0$ such that
\BQN\label{holderf}
\max_{i=1,\dots, n}\sup_{s\neq t, s,t\in [0,T]}\frac{|h_i(t)-h_i(s)|}{|t-s|^{\mu_2}}<\IF.
\EQN
\BT\label{Thm2}
Suppose that  ${\bf X}(t),t\geq 0$ is a centered vector-valued Gaussian process with independent coordinates satisfying (\ref{eq:sigt0})-(\ref{holderx}),  and  $\vk{h}(t), \ t\geq 0$ is a  continuous vector function over $[0,T]$ satisfying (\ref{eq:gtt0})-(\ref{holderf}).
Then
\BQNY
\pk{\exists_{t\in[0,T]}\LT(\vk{X}(t)+{\vk{h}}(t)\RT) >u\vk{1}}&\sim& u^{(\frac{2}{\alpha}-\frac{2}{\beta})_{+}}\prod_{i=1}^n\Psi\LT(\frac{u-h_i(t_0)}
{\sigma_i(t_0)}\RT)
\\
&&\times\LT\{
\begin{array}{ll}
\mathcal{H}_{\alpha,\frac{\vk{a}}{\vk{\sigma}^2(t_0)}
\mathbf{I}_{\{\vk{\alpha}=\alpha\vk{1}\}}}\int_{\cLa{q}}^\IF e^{-\sum_{i=1}^n f_i(x)}dx,&\ \ \text{if}\ \alpha<\beta,\\
\mathcal{P}^{\vk{f}}_{\alpha,\frac{\vk{a}}{\vk{\sigma}^2(t_0)}
\mathbf{I}_{\{\vk{\alpha}=\alpha\vk{1}\}}}[\cLa{q},\IF), &\ \ \text{if}\ \alpha=\beta,\\
1,&\ \ \text{if}\ \alpha>\beta,
\end{array}
\RT.
\EQNY
where $\alpha=\min_{1\leq i\leq n}\alpha_i$, $\beta=\min_{1\leq i\leq n}\min(\beta_i, 2\gamma_i\II_{\{c_i\neq0\}}+\IF\II_{\{c_i=0\}})$, $\vk{a}=(a_1,\dots, a_n)$, $\vk{\sigma}(t_0)=(\sigma_1(t_0),\dots, \sigma_n(t_0))$, $\vk{f}=(f_1,\dots, f_n)$ with
$f_i(t)=\frac{b_i}{\sigma_i^3(t_0)}|t|^{\beta_i} \II_{\{\beta_i=\beta\}}+\frac{c_i}{\sigma_i^2(t_0)}|t|^{\gamma_i} \mathbb{I}_{\{2\gamma_i=\beta\}}$, and
\begin{align}\label{II}
\cLa{q}= \LT\{
\begin{array}{ll}
-\IF,& \hbox{if} \  t_0\in(0,T),\\
0,& \hbox{if} \  t_0=0\ \hbox{or}\ t_0=T.
\end{array}
\RT.
\end{align}
\ET

\begin{remark}
 If $n=1$ and $h_1(t)\equiv 0$,  then  \netheo{Thm2} covers the classical
 Piterbarg-Prisja{\v{z}}njuk result; see \cite{PP78}.
\end{remark}

In the following corollary we apply  \netheo{Thm2} for
the analysis of exact asymptotics of
$\tau_u=\inf\{t\geq 0: \LT(\vk{X}(t)+\vk{h}(t)\RT) >u\vk{1}\}$, as $u\to\infty$,
conditioned that $\tau_u\le T$.

\BK\label{Korr1}
Under the same assumptions  as in \netheo{Thm2} with $t_0=T$, we have for $x\in(0,\IF)$, as $u\rw\IF$,
\BQN
\pk{(T-\tau_u)u^{2/\beta}\leq x\Big| \tau_u\leq T}\sim
\LT\{
\begin{array}{ll}
\int_{0}^x e^{-\sum_{i=1}^n f_i(t)}dt\Big/\int_{0}^\IF e^{-\sum_{i=1}^n f_i(t)}dt,&\ \ \text{if}\ \alpha<\beta,\\
\mathcal{P}^{\vk{f}}_{\alpha,\frac{\vk{a}}{\vk{\sigma}^2(t_0)}
\mathbf{I}_{\{\vk{\alpha}=\alpha\vk{1}\}}}[0,x]\Big/\mathcal{P}^{\vk{f}}_{\alpha,\frac{\vk{a}}{\vk{\sigma}^2(t_0)}
\mathbf{I}_{\{\vk{\alpha}=\alpha\vk{1}\}}}[0,\IF), &\ \ \text{if}\ \alpha=\beta,\\
1,&\ \ \text{if}\ \alpha>\beta.
\end{array}
\RT.
\EQN
\EK

We give a short proof of Corollary \ref{Korr1} in Appendix.

\def\ubi{\underline{b}_i}
\def\obi{\overline{b}_i}
\COM{Next, let us focus on the
case of purely centered vector-valued Gaussian processes, i.e.,  $\vk{h}(t)\equiv\vk{0}$.
The following theorem extends findings of
\cite{Tabis}, where only the case of
$\beta_i=\beta, 1\leq i\leq n$ in (\ref{bi1}) was  considered.
%\BT\label{Thmm3}
\begin{example}
Let ${\bf X}(t),t\geq 0$ be a centered vector-valued Gaussian process with independent coordinates and continuous trajectories. Assume that $g(t)=\sum_{i=1}^n\frac{1}{\sigma_i^2(t)}$, with $\sigma_i^2(t)$ representing the variance of $i-$th coordinates,  attains its maximum  over $[0, T]$ at the unique point $t_0\in[0,T]$ satisfying
\BQN\label{bi1}
\sigma_i(t)=\sigma_i(t_0)-\ubi|t-t_0|^{\beta_i}\mathbb{I}_{\{t\leq 0\}}-\obi|t-t_0|^{\beta_i}\mathbb{I}_{\{t> 0\}}+o(|t-t_0|^{\beta_i}),\ \
\ \ t\rw t_0
\EQN
for some constants $\beta_i>0,\underline{b}_i,\overline{b}_i\in\R $.  Moreover, we further assume that   \eqref{eq:rst} and (\ref{holderx}) are satisfied.
If further
$$\sum_{i=1}^n\frac{\obi\II_{\{\beta_i=\beta\}}}{\sigma_i^3(t_0)}>0, \quad \sum_{i=1}^n\frac{\ubi\II_{\{\beta_i=\beta\}}}{\sigma_i^3(t_0)}>0,$$
 then by Theorem \ref{Thm2}
\BQNY
\pk{\exists_{t\in[0,T]}\vk{X}(t) >u\vk{1}}\sim u^{(\frac{2}{\alpha}-\frac{2}{\beta})_{+}}\prod_{i=1}^n\Psi\LT(\frac{u}{\sigma_i(t_0)}\RT)
\LT\{
\begin{array}{ll}
\mathcal{H}_{\alpha,\frac{\vk{a}}{\vk{\sigma}^2(t_0)}
\mathbf{I}_{\{\vk{\alpha}=\alpha\vk{1}\}}}\int_{Q}^\IF e^{-\sum_{i=1}^n f_i(x)}dx,&\ \ \text{if}\ \alpha<\beta,\\
\mathcal{P}^{\vk{f}}_{\alpha,\frac{\vk{a}}{\vk{\sigma}^2(t_0)}
\mathbf{I}_{\{\vk{\alpha}=\alpha\vk{1}\}}}[Q,\IF), &\ \ \text{if}\ \alpha=\beta,\\
1,&\ \ \text{if}\ \alpha>\beta,
\end{array}
\RT.
\EQNY
where $\alpha=\min_{1\leq i\leq n}\alpha_i$, $\beta=\min_{1\leq i\leq n}\beta_i$, $f_i(t)=\frac{1}{\sigma_i^3(t_0)}\LT(\ubi|t|^{\beta_i}\mathbb{I}_{\{t\leq 0\}}+\obi|t|^{\beta_i}\mathbb{I}_{\{t\leq 0\}}\RT)\II_{\{\beta_i=\beta\}}$, and $Q$ is defined in \eqref{II}.
\end{example}}

%%%%%%%%%%%%%%%%%%%%%%%%%%%%%%%%%%%%%%%%%%%%%55

\subsection{Locally-stationary coordinates}

\cLa{Suppose that  for each $i=1,...,n$, $X_i$ is a
centered locally-stationary Gaussian process} with continuous trajectories, that is process with
unit  variance and correlation function $r_{i}(\cdot,\cdot), 1\leq i\leq n$ satisfying
 \BQN \label{stationaryR0}
r_i(t,t+s) = 1 - a_i(t) \abs{s}^{\alpha_i} + o(\abs{s}^{\alpha_i}), \  \ s \to 0
\EQN
uniformly with respect to $t\in[0,T]$,
where  $\alpha_i \in (0,2]$, and
$a_i(t)\in (0,\infty)$ is a  positive continuous function  on $[0,T]$.
Further, we suppose that
 \BQN \label{stationaryR2}
  r_i(s,t)< 1, \  \forall s,t\in[0,T]\ \mathrm{and\ } s \not=t.
\EQN

We refer to e.g., \cite{Ber74, Berman92, Hus90, Pit96}  for the investigation of extremes of one-dimensional  locally-stationary  Gaussian processes
under the above conditions.

Denote by
$$H=\bigcap_{i=1}^n\LT\{s\in[0,T]:h_i(s)=h_{m,i}:=\max_{t\in[0,T]}h_i(t)\RT\}.
$$

\BT\label{Thm3}
Let ${\bf X}(t),t\in[0,T]$  be a locally stationary vector-valued Gaussian process satisfying (\ref{stationaryR0}) and (\ref{stationaryR2}).
 Moreover, assume that $\vk{h}(t)$ is a vector function satisfying (\ref{holderf}) and $\alpha=\min_{1\leq i\leq n}\alpha_i.$\\
i) If $H=\{t_0\}$ and \eqref{eq:gtt0} holds with $c_i\geq 0$ and $\max_{1\leq i\leq n}c_i>0$, then
\BQNY
\pk{\exists_{t\in[0,T]}\LT(\vk{X}(t)+{\vk{h}}(t)\RT) >u\vk{1}}
\sim
u^{(\frac{2}{\alpha}-\frac{1}{\gamma})_{+}}\prod_{i=1}^n\Psi\LT(u-h_{m,i}\RT)
\LT\{
\begin{array}{ll}
\mathcal{H}_{\alpha,\vk{a}(t_0)\mathbf{I}_{\{\vk{\alpha}=\alpha\vk{1}\}}}\int_{q}^\IF e^{-\sum_{i=1}^n f_i(x)}dx,&\ \ \text{if}\ \alpha<2\gamma,\\
\mathcal{P}^{\vk{f}}_{\alpha,\vk{a}(t_0)\mathbf{I}_{\{\vk{\alpha}=\alpha\vk{1}\}}}[q,\IF), &\ \ \text{if}\ \alpha=2\gamma,\\
1,&\ \ \text{if}\ \alpha>2\gamma,
\end{array}
\RT.
\EQNY
where  $\gamma=\min_{1\leq i\leq n}( \gamma_i\II_{\{c_i\neq0\}}+\IF\II_{\{c_i=0\}})$,
$f_i(t)=c_i|t|^\gamma \mathbb{I}_{\{\gamma_i=\gamma\}}$, and $q$ is given by \eqref{II}.\\
ii) If $H=[A,B]\subset[0,T]$ with $A>B$,
 then
\BQNY
\pk{\exists_{t\in[0,T]}\LT(\vk{X}(t)+{\vk{h}}(t)\RT) >u\vk{1}}
\sim
\int_{A}^B\mathcal{H}_{\alpha,\vk{a}(t)\mathbf{I}_{\{\vk{\alpha}=\alpha\vk{1}\}}}dt\
u^{\frac{2}{\alpha}}
\prod_{i=1}^n\Psi\LT(u-h_{m,i}\RT).
\EQNY
\COM{(iii) If $H=\{t_1,\ldots, t_d\}\subset(0,T)$ with $t_1<t_2<\ldots<t_d$, $d\geq1$, and
\BQNY
h_i(t)=h_{m,i}-c_{i,j}|t-t_j|^{\gamma_{i,j}}(1+o(1)),\  t\rw t_j, j=1,\ldots d,\ i=1,\ldots,n,
\EQNY
for some constants $\gamma_{i,j}>0$, $c_{i,j}\ge 0$ and $\max_{1\le i\le n}c_{i,j}>0, j=1,\ldots,d$. Then
\BQNY
\pk{\exists_{t\in[0,T]}\LT(\vk{X}(t)+{\vk{h}}(t)\RT)>u\vk{1}}\sim
\LT\{
\begin{array}{ll}
\sum_{j=1}^d \II_{\{\gamma^j=\gamma\}}u^{\frac{2}{\alpha}-\frac{1}{\gamma}}
\mathcal{H}_{\alpha,\vk{a}(t_j)\mathbf{I}_{\{\vk{\alpha}=\alpha\vk{1}\}}}
\int_{-\IF}^{\IF}e^{-\sum_{i=1}^nf_i^j(x)}dx\prod_{i=1}^n\Psi\LT(u-h_{m,i}\RT),& 2\gamma>\alpha,\\
\sum_{j=1}^d \LT(\II_{\{2\gamma^j<\alpha\}}+\II_{\{2\gamma^j=\alpha\}}
\mathcal{P}_{\alpha,\vk{a}(t_j)\mathbf{I}_{\{\vk{\alpha}=\alpha\vk{1}\}}}^{\vk{f}^j}(-\IF,\IF)\RT)
\prod_{i=1}^n\Psi\LT(u-h_{m,i}\RT),& 2\gamma\le \alpha,
\end{array}
\RT.
\EQNY
where $\alpha=\min_{1\leq i\leq n}\alpha_i$, $\gamma^j=\min_{1\le i\le n}\LT(\gamma_{i,j}\mathbb{I}_{\{c_{i,j}\neq 0\}}+\IF\mathbb{I}_{\{c_{i,j}=0\}}\RT)$,$\gamma=\max_{1\leq j\leq d}\gamma^j$,
$f^j_i(t)=c_{i,j}|t|^{\gamma^j} \mathbb{I}_{\{\gamma_{i,j}=\gamma^j\}}$, and $\vk{f}^j=(f^j_1,\ldots,f^j_n)$.}
\ET
\COM{\begin{example}
Let ${\bf X}(t),t\in[0,T]$ be a centered vector-valued Gaussian process
with independent coordinates.
Suppose that $X_i$, $i=1,...,n$ are  stationary Gaussian processes
with continuous trajectories,
unit  variance  and correlation functions
$r_{i}(\cdot,\cdot), 1\leq i\leq n$ satisfying
\eqref{stationaryR0} and \eqref{stationaryR2} with $a_i(t)\equiv a_i$.
If
$c_i\geq 0$ and $\max_{1\leq i\leq n}c_i>0$, then we have
\BQNY
\pk{\exists_{t\in[0,T]}\LT(\vk{X}(t)+{\vk{c}}t\RT) >u\vk{1}}\sim u^{(\frac{2}{\alpha}-1)_{+}}\prod_{i=1}^n\Psi\LT(u-c_i T\RT)
\LT\{
\begin{array}{ll}
\frac{1}{\sum_{i=1}^n c_i}\mathcal{H}_{\alpha,\vk{a}\mathbf{I}_{\{\vk{\alpha}=\alpha\vk{1}\}}},&\ \ \text{if}\ \alpha<2,\\
\mathcal{P}^{\vk{f}}_{\alpha,\vk{a}\mathbf{I}_{\{\vk{\alpha}=\alpha\vk{1}\}}}[0,\IF), &\ \ \text{if}\ \alpha=2
\end{array}
\RT.
\EQNY
where $\alpha=\min_{1\leq i\leq n}\alpha_i$,
$f_i(t)=c_i|t| $.
The proof of the above formula follows from (i) of \netheo{Thm3},
where we use that
$H=\bigcap_{i=1}^n\LT\{s\in[0,T]: c_i s=\max_{t\in[0,T]}c_i t\RT\}=\{T\}.$
\end{example}}

\COM{
\BT\label{Thm4}
Let $\{{\bf X}(t),t\in[0,T]\}$ be given as in \netheo{Thm3} and $h(t)$ be bound measurable functions. Further, $h(t)$ attains its maximum $h_m=\max_{t\in[0,T]}h(t)$ at $t\in E$.\\
(1) If $E=\{t_1,\ldots, t_d\}\subset(0,T)$ with $t_1<t_2<\ldots<t_d$, $d\geq1$, and
\BQNY
h(t)=h_m-c_j|t-t_j|^{\gamma_j}(1+o(1)),\  t\rw t_j, j=1,\ldots d,
\EQNY
for some positive constants $\gamma_j$ and $c_j$, then
\BQNY
\pk{\sup_{t\in[0,T]}\LT(\min_{1\le i\le n}X_i(t)\RT)+h(t)>u}\sim u^{(\frac{2}{\alpha}-\frac{1}{\gamma})_{+}}\prod_{i=1}^n\Psi\LT(u-h_i(t_0)\RT)
\LT\{
\begin{array}{ll}
\mathcal{H}_{\alpha,\vk{a}\mathbf{I}_{\{\vk{\alpha}=\alpha\vk{1}\}}}\int_{-\hat{I}}^\IF e^{-\sum_{i=1}^n f_i(x)}dx,&\ \ \text{if}\ \alpha<2\gamma,\\
\mathcal{P}^{\vk{f}}_{\alpha,\vk{a}\mathbf{I}_{\{\vk{\alpha}=\alpha\vk{1}\}}}[-\hat{I},\IF), &\ \ \text{if}\ \alpha=2\gamma,\\
1,&\ \ \text{if}\ \alpha>2\gamma,
\end{array}
\RT.
\EQNY
where $\alpha=\min_{1\leq i\leq n}\alpha_i$, $\gamma=\max_{1\leq j\leq d}\gamma_j$,
$f_i(t)=c_i|t|^\gamma \mathbb{I}_{\{\gamma_i=\gamma\}}$, and $ \hat{I}= \IF \II_{\{t_0\in(0,T)\}}$.\\
(2)If $E=[A,B]\subset[0,T]$ with $A>B$,
when $A>0$,
\BQNY
h(t)=h_m-c_{A}|t-A|^{\gamma_{A}}(1+o(1)),\  t\uparrow A,
\EQNY
and when $B<T$,
\BQNY
h(t)=h_m-c_{B}|t-B|^{\gamma_{B}}(1+o(1)),\  t\downarrow B,
\EQNY
for some positive constants $c_A,c_B,\gamma_A,\gamma_B$, then
\BQNY
\pk{\sup_{t\in[0,T]}\LT(\min_{1\le i\le n}X_i(t)\RT)+h(t)>u}\sim u^{\frac{2}{\alpha}}\prod_{i=1}^n\Psi\LT(u-h_m\RT)
\LT\{
\begin{array}{ll}
\mathcal{H}_{\alpha,\vk{a}\mathbf{I}_{\{\vk{\alpha}=\alpha\vk{1}\}}}\int_{-\hat{I}}^\IF e^{-\sum_{i=1}^n f_i(x)}dx,&\ \ \text{if}\ \alpha<2\gamma,\\
\mathcal{P}^{\vk{f}}_{\alpha,\vk{a}\mathbf{I}_{\{\vk{\alpha}=\alpha\vk{1}\}}}[-\hat{I},\IF], &\ \ \text{if}\ \alpha=2\gamma,\\
1,&\ \ \text{if}\ \alpha>2\gamma,
\end{array}
\RT.
\EQNY
where $\alpha=\min_{1\leq i\leq n}\alpha_i$, $\gamma=\min_{1\leq i\leq n}( \gamma_i\II_{\{c_i\neq0\}}+\IF\II_{\{c_i=0\}})$,
$f_i(t)=c_i|t|^\gamma \mathbb{I}_{\{\gamma_i=\gamma\}}$, and $ \hat{I}= \IF \II_{\{t_0\in(0,T)\}}$.
\ET}
Similarly to Corollary  \ref{Korr1}, we get the asymptotics of $\tau_u$ for locally-stationary coordinates of $\vk{X}$.

\BK\label{Korr2}
Under the same assumptions as in i) of \netheo{Thm3}, with $t_0=T$, we have for $x\in(0,\IF)$, as $u\rw\IF$,
\BQN
\pk{(T-\tau_u)u^{1/\gamma}\leq x\Big| \tau_u\leq T}\sim
\LT\{
\begin{array}{ll}
\int_{0}^x e^{-\sum_{i=1}^n f_i(t)}dt\Big/\int_{0}^\IF e^{-\sum_{i=1}^n f_i(t)}dt,&\ \ \text{if}\ \alpha<2\gamma,\\
\mathcal{P}^{\vk{f}}_{\alpha,\vk{a}(t_0)
\mathbf{I}_{\{\vk{\alpha}=\alpha\vk{1}\}}}[0,x]
\Big/\mathcal{P}^{\vk{f}}_{\alpha,\vk{a}(t_0)
\mathbf{I}_{\{\vk{\alpha}=\alpha\vk{1}\}}}[0,\IF), &\ \ \text{if}\ \alpha=2\gamma,\\
1,&\ \ \text{if}\ \alpha>2\gamma.
\end{array}
\RT.
\EQN
\EK

\subsection{A simultaneous ruin model}
Consider  portfolio $\vk{U}(t)=(U_1(t)\ldot U_n(t))$, where
$$\vk{U}(t)=u\vk{d}+\vk{c}t-\vk{B}_{\vk{\alpha}}(t), \quad t\ge0,$$
with $\vk{c}=(c_1,\cdots,c_n)\in \R^n$, $\vk{d}=(d_1,\cdots,d_n)>\vk{0}$ and $B_{\alpha_i}(t),\ 1\leq i\leq n$,
independent standard fractional Brownian motions with  variance
$\Var(B_{\alpha_i}(t))= t^{\alpha_i}$ for  $\alpha_i\in(0,2], \ 1\leq i\leq n$, respectively.
The corresponding simultaneous ruin probability over $[0,T]$ is defined as
$$ \pk{\exists_{t\in[0,T]}\vk{U}(t) < \vk{0}} $$
and  the simultaneous ruin time   $\tau_u:=\inf\{t\geq0:\vk{U}(t) < \vk{0}\}$.
We refer to, e.g., \cite{Mic98} for theoretical justification of the use of fractional Brownian motion
as the approximation of {\it the claim} process in risk theory.

In the following proposition we present exact
asymptotics of the simultaneous ruin probability and the conditional simultaneous ruin time $\tau_u|\tau_u<T$, as $u\to \IF$.

\begin{prop}\label{prop1}
For $T\in (0,\infty)$, $\alpha=\min_{1\leq i\leq n}\alpha_i$,
$b_i=\frac{d_i^2}{2T^{2\alpha_i}}$ and $f_i(t)=\frac{\alpha_id_i^2}{2T^{\alpha_i+1}}t$,
as $u\rw\IF$, we have
\begin{align}\label{e.1}
\pk{\exists_{t\in[0,T]}\vk{U}(t) < \vk{0}}\sim &u^{(\frac{2}{\alpha}-2)_{+}}\prod_{i=1}^n\Psi\LT(\frac{d_i u+c_iT}{T^{{\alpha_i}/2}}\RT)\\
&\times\LT\{
\begin{array}{ll}
\LT(\sum_{i=1}^n\frac{\alpha_id_i^2}{2T^{\alpha_i+1}}\RT)^{-1}
\mathcal{H}_{\alpha, \cLa{\vk{b}}\mathbf{I}_{\{\vk{\alpha}=\alpha\vk{1}\}}},
&\ \ \text{if}\ \alpha<1,\\
\mathcal{P}^{\vk{f}}_{\alpha,\cLa{\vk{b}}
\mathbf{I}_{\{\vk{\alpha}=\alpha\vk{1}\}}}[0,\IF), &\ \ \text{if}\ \alpha=1,\\
1,&\ \ \text{if}\ \alpha>1
\end{array}
\RT.\nonumber
\end{align}
and
for $x\in (0,\IF)$
\BQN
\pk{(T-\tau_u)u^2\leq x\Big| \tau_u\leq T}\sim
\LT\{
\begin{array}{ll}
1-e^{-\LT(\sum_{i=1}^n\frac{\alpha_id_i^2}{2T^{\alpha_i+1}}\RT)x},&\ \ \text{if}\ \alpha<1,\\
{\mathcal{P}^{\vk{f}}_{\alpha,\vk{b}
\mathbf{I}_{\{\vk{\alpha}=\alpha\vk{1}\}}}[0,x]}\Big/
{\mathcal{P}^{\vk{f}}_{\alpha,\vk{b}
\mathbf{I}_{\{\vk{\alpha}=\alpha\vk{1}\}}}[0,\IF)}, &\ \ \text{if}\ \alpha=1,\\
1,&\ \ \text{if}\ \alpha>1.
\end{array}
\RT.
\label{e.2}
\EQN
\end{prop}

Specifically, Proposition \ref{prop1}
allows us to get exact asymptotics for
multidimensional counterpart of
the classical Brownian risk model \cite{Igl69}.
For simplicity  we focus on 2-dimensional case.
Let $\vk{B}(t):=(B^{(1)}(t),B^{(2)}(t) )$,
where $B^{(1)}(t)$ and $B^{(2)}(t)$ are two independent standard Brownian motions, $\vk{c}=(c_1,c_2)\in \R^2$
and $\vk{d}=(d_1,d_2)\in \R^2_+$.
Then
we have, as $u\to\infty$,
\begin{align*}
\pk{\exists_{t\in[0,T]}
\left(\begin{array}{l}
d_1u+c_1t-B^{(1)}(t)  \\
d_2u+c_2t-B^{(2)}(t)
\end{array}
\right)
\le
\left(\begin{array}{l}
0  \\
0
\end{array}
\right)
}\sim \mathcal{P}^{\vk{b}t}
_{1,\vk{b}}[0,\IF)
\Psi\LT(\frac{d_1u+c_1T}{T^{1/2}}\RT)\Psi\LT(\frac{d_2u+c_2T}{T^{1/2}}\RT)
\end{align*}
and for $x\in(0,\IF)$
\begin{align*}
\pk{(T-\tau_u)u^2\leq x\Big| \tau_u\leq T}\sim \mathcal{P}^{\vk{b}t}
_{1,\vk{b}}[0,x]\Big/\mathcal{P}^{\vk{b}t}
_{1,\vk{b}}[0,\IF),
\end{align*}
where $\vk{b}=\LT(\frac{d^2_1}{2T^2},\frac{d_2^2}{2T^2}\RT)$.

\section{Proofs}
 Before proceeding to the proof of \netheo{PreThm1},  we  present two  lemmas which play an important role in the proof of \netheo{PreThm1}.
 The first one is a vector-valued version of the uniform Pickands-Piterbarg lemma while the second one gives an upper bound for  the double maximum of vector-valued Gaussian process. Hereafter,  we denote by $\mathbb{C}_l, l\in\N$  some positive constants that may differ from line to line. Moreover, the notation $f(u, S, \epsilon)\sim g(u)$ as $u\rw\IF, S\rw\IF, \epsilon\rw 0$,  means that  $\lim_{\epsilon\rw 0}\lim_{S\rw\IF}\lim_{u\rw\IF}\frac{f(u,S,\epsilon)}{g(u)}=1$. \\
For $\vk{b}\geq\vk{0}$, $\lambda_{i}\in[0,\IF)$,  and $-\IF<S_1<S_2<\IF$, define a vector-valued Gaussian process $\vk{Z}_u(t)=(Z_{u,1}(t),\ldots,Z_{u,n}(t))$ by
\BQN\label{zu}
Z_{u,i}(t)=\frac{\xi_i(t)}{1+b_iu^{-2}f_{i}(u^{\lambda_{i}}t)},\ t\in[S_1,S_2], \ i=1,\ldots, n,
\EQN
where $\vk{\xi}(t)=(\xi_1(t), \dots, \xi_n(t)), t\in\R$ is a vector-valued Gaussian process with independent stationary coordinates,  continuous sample paths, unit variance and correlation function $r_i(\cdot)$  on $i$-th coordinate, $1\leq i\leq n$,  satisfying
\BQN\label{rp}
1-r_i(t)=a_i\abs{t}^{\alpha_i}\bl{(1+o(1))},
\EQN
for $a_i>0$ and $\alpha_i\in(0,2]$, and $f_{i}(t), 1\leq i\leq n$ are some continuous functions.  We suppose that the threshold vector \cLa{$\vk{m}_u(k)=(m_{u,1}(k)\ldot m_{u,n}(k))$} satisfies
\BQN\label{Mku}
\lim_{u\to\IF}\sup_{k\in K_u}\LT|\frac{1}{u}\vk{m}_u(k)-\vk{c}\RT|=0,\ \ \vk{c}>\vk{0},
\EQN with $K_u$  an index set.

Denote by $$\alpha=\min_{1\leq i\leq n}\alpha_i, \quad  \lambda=\max_{1\leq i\leq n}(\lambda_{i}\II_{\{b_i\neq 0\}})>0, \quad
\widetilde{\vk{f}}(t)=(\widetilde{f}_1(t),\ldots,\widetilde{f}_n(t)), \quad \text{with}\quad
\widetilde{f}_i\LT(t\RT)=f_i\LT(t\RT)\mathbb{I}_{\{\lambda_i=\lambda\}}.$$

\BEL\label{pan2} Let $\vk{Z}_u(t)$ be defined in (\ref{zu}) and $\vk{m}_u(k)$ satisfy (\ref{Mku}).\\
i) If $\lambda\leq 2/\alpha$, then
\BQNY
\lim_{u\rw\IF}\sup_{k\in K_u}\LT| %\LT(\sqrt{2\pi}M_{u,i}(k)e^{M_{u,i}^2(k)/2}\RT)
\frac{\pk{\exists_{t\in
[u^{-2/\alpha}S_1,u^{-2/\alpha}S_2]}\vk{Z}_u(t)>\vk{m}_u(k)}}{\prod_{i=1}^n
\Psi(M_{u,i}(k))}-\mathcal{R}_{\lambda}^f[S_1,S_2]\RT|=0,
\EQNY
where
\BQNY
\mathcal{R}_{\lambda}^f[S_1,S_2]= \LT\{
\begin{array}{ll}
\mathcal{P}_{\alpha,\vk{a}{\vk{c}}^2\mathbf{I}_{\{\vk{\alpha}=\alpha\vk{1}}\}}^{\vk{c}^2
\widetilde{\vk{f}}}[S_1,S_2],&\ \hbox{if}\ \lambda=2/\alpha,\\
\mathcal{P}_{\alpha,\vk{a}{\vk{c}}^2\mathbf{I}_{\{\vk{\alpha}=\alpha\vk{1}}\}}^{\vk{c}^2
\widetilde{\vk{f}}(0)}[S_1,S_2],&\ \hbox{if}\ \lambda<2/\alpha,\\
\mathcal{H}_{\alpha,\vk{a}{\vk{c}}^2\mathbf{I}_{\{\vk{\alpha}=\alpha\vk{1}}\}}[S_1,S_2],&\ \hbox{if}\  \vk{b}=\vk{0}.
\end{array}
\RT.
\EQNY
ii) If $\lambda>2/\alpha$, then
\BQNY
\lim_{u\rw\IF}\sup_{k\in K_u}\LT| %\LT(\sqrt{2\pi}M_{u,i}(k)e^{M_{u,i}^2(k)/2}\RT)
\frac{\pk{\exists_{t\in
[u^{-\lambda}S_1,u^{-\lambda}S_2]}\vk{Z}_u(t)>\vk{m}_u(k)}}{\prod_{i=1}^n
\Psi(m_{u,i}(k))}-
\mathcal{P}_{\alpha,\vk{0}}^{\vk{c}^2\widetilde{\vk{f}}}[S_1,S_2]\RT|=0.
\EQNY
\COM{where
\BQNY
\mathcal{R}_{\lambda}^f[S_1,S_2]=\int_{\R^n}e^{\sum_{i=1}^n w_i}\mathbb{I}_{\{\exists_{t\in[S_1,S_2]}
-\vk{c}^2\widehat{\vk{f}}(t)>\vk{w}\}}d\vk{w},
%= \sup_{t\in[S_1,S_2]}e^{-\sum_{i=1}^nc_i^2f_i(t)},
\EQNY
and  $\mathcal{R}_{\lambda}^f[S_1,S_2]=1$ if $\vk{b}=\vk{0}$.}
\EEL

{\bf Proof}.
i) Suppose that $\lambda\leq 2/\alpha$.
Conditioning on \cLa{$\LT\{\vk{\xi}(0)=
\vk{m}_u(k)-\frac{\vk{w}}{\vk{m}_u(k)}\RT\}, \vk{w}\in {\R}^n$}, we have for all $u$ large enough
\BQNY
&&\frac{\pk{\exists_{t\in
[u^{-2/\alpha}S_1,u^{-2/\alpha}S_2]}\vk{Z}_u(t)>\vk{m}_u(k)}}{\prod_{i=1}^n
\Psi(m_{u,i}(k))}\\
&&=\frac{1}{\prod_{i=1}^n\sqrt{2\pi}m_{u,i}(k)\Psi(m_{u,i}(k))}\int_{{\R}^n}e^{-\frac{1}{2}
\sum_{i=1}^n\LT(m_{u,i}(k)-\frac{w_i}{m_{u,i}(k)}\RT)^2}\\
&&\quad\times\pk{\exists_{t\in
[S_1,S_2]}\vk{Z}_u(u^{-2/\alpha}t)>\vk{m}_u(k)\Bigl\lvert \vk{\xi}(0)=
\vk{m}_u(k)-\frac{\vk{w}}{\vk{m}_u(k)} }d\vk{w}\\
&&=\LT(\prod_{i=1}^n\frac{e^{-\frac{\left(m_{u,i}(k)\right)^2}{2}}}{\sqrt{2\pi}m_{u,i}(k)
\Psi(m_{u,i}(k))}\RT)\int_{{\R}^n}
e^{\sum_{i=1}^n\left(w_i-\frac{w_i^2}{2\left(m_{u,i}(k)\right)^2}\right)}\pk{\exists_{t\in
[S_1,S_2]}\vk{\mathcal{X}}^{\vk{w}}_u(t,k)>\vk{w} }d\vk{w}\\
&&=\LT(\prod_{i=1}^n\frac{e^{-\frac{\left(m_{u,i}(k)\right)^2}{2}}}{\sqrt{2\pi}m_{u,i}(k)
\Psi(m_{u,i}(k))}\RT)I_{u,k},
\EQNY
where $\vk{\mathcal{X}}^{\vk{w}}_u(t,k)=(\mathcal{X}^{w}_{u,1}(t,k),\ldots,
\mathcal{X}^{w}_{u,n}(t,k))$ with
\BQNY
\mathcal{X}^{w}_{u,i}(t,k)=m_{u,i}(k)(Z_{u,i}(u^{-2/\alpha}t)-m_{u,i}(k))+w_i\Bigl
\lvert\xi_i(0)=
m_{u,i}(k)-\frac{w_i}{m_{u,i}(k)}.
\EQNY
By (\ref{Mku}), it follows that
$$\lim_{u\rw\IF}\sup_{k\in K_u}\left|\LT(\prod_{i=1}^n\frac{e^{-\frac{\left(m_{u,i}(k)\right)^2}{2}}}{\sqrt{2\pi}
m_{u,i}(k)\Psi(m_{u,i}(k))}\RT)-1\right|=0.$$
Thus in order to establish the proof, it suffices to prove that
\BQN\label{IR}
\lim_{u\rw\IF}\sup_{k\in K_u}\LT|I_{u,k}-\mathcal{R}_{\lambda}^f[S_1,S_2]\RT|=0.
\EQN
It follows that, for each $W>0$, with $\widetilde{W}^n=[-W,W]^n$  and $\widetilde{W}_j^n=\{\vk{w}\in\R^n\big|w_j\in(-\IF,-W)\cup(W,\IF)\}$,
\BQNY
&&\sup_{k\in K_u}\LT|I_{u,k}-\mathcal{R}_{\lambda}^f[S_1,S_2]\RT|\nonumber\\
&&\leq\sup_{k\in K_u}\LT|\int_{\widetilde{W}^n}\LT[e^{\sum_{i=1}^n\left(w_i-\frac{w_i^2}{2m_{u,i}^2(k)}\right)}
\pk{\exists_{t\in
[S_1,S_2]}\vk{\mathcal{X}}^{\vk{w}}_u(t,k)>\vk{w} }
-e^{\sum_{i=1}^n w_i}\pk{\exists_{t\in[S_1,S_2]}\vk{\zeta}(t)>\vk{w}}\RT]d\vk{w}\RT|\nonumber\\
&&\quad+\sum_{j=1}^n\sup_{k\in K_u}\int_{\widetilde{W}_j^n}e^{\sum_{i=1}^n w_i}\pk{\exists_{t\in
[S_1,S_2]}\vk{\mathcal{X}}^{\vk{w}}_u(t,k)>\vk{w} }d\vk{w}\label{pi}\\
&&\quad+\sum_{j=1}^n\int_{\widetilde{W}_j^n}e^{\sum_{i=1}^n w_i}\pk{\exists_{t\in[S_1,S_2]}\vk{\zeta}(t)>\vk{w}}d\vk{w}\\
&&:=I_1(u)+I_2(u)+I_3(u),
\EQNY
where  $\vk{\zeta}(t)=(\vk{c}\sqrt{2\vk{a}}\vk{B}_{\alpha}-\vk{a}{\vk{c}}^2|t|^\alpha)
\mathbf{I}_{\{\vk{\alpha}=\alpha\vk{1}\}}
-\vk{c}^2\widetilde{\vk{f}}(t\II_{\{\lambda=2/\alpha\}})$.

Next, we give upper bounds for $I_i(u), i=1,2,3$. We begin with the weak convergence of process
$\vk{\mathcal{X}}^{\vk{w}}_u(t,k)$.

\underline{Weak convergence of $\vk{\mathcal{X}}^{\vk{w}}_u(t,k)$}.
Direct calculation shows that
\BQNY\label{XM11}
\EF{(1+b_iu^{-2}f_{i}(u^{\lambda_{i}}t))\mathcal{X}^w_{u,i}(t,k)}
&=&-m_{u,i}^2(k)\LT(1-r_i(u^{-2/\alpha}t)+b_iu^{-2}f_{i}(u^{\lambda_{i}-2/\alpha}t)\RT)\\
&&+w_i\LT(1-r_i(u^{-2/\alpha}t)+b_iu^{-2}f_{i}(u^{\lambda_{i}-2/\alpha}t)\RT),
\EQNY
and
\BQNY
&&\Var\LT((1+b_iu^{-2}f_{i}(u^{\lambda_{i}}t))\mathcal{X}^w_{u,i}(t,k)-(1+b_iu^{-2}
f_{i}(u^{\lambda_{i}}t'))\mathcal{X}^w_{u,i}(t',k)\RT)\\
&&=m_{u,i}^2(k)
\LT(Var\LT(\xi_i(u^{-2/\alpha}t)-\xi_i(u^{-2/\alpha}t')\RT)-
\LT(r_i(u^{-2/\alpha}t)
-r_i(u^{-2/\alpha}t')\RT)^2\RT).
\EQNY
By (\ref{rp}) and (\ref{Mku}), it follows  that
\BQN\label{XM12}
\EF{(1+b_iu^{-2}f_{i}(u^{\lambda_{i}}t))\mathcal{X}^w_{u,i}(t,k)}\rw-c_i^2a_i|t|^\alpha
\II_{\{\alpha_i=\alpha\}}
-c_i^2\LT(\widetilde{f}_{i}(t\II_{\{\lambda=2/\alpha\}})\RT),
\EQN
 as $u\rw\IF$, uniformly with respect to $t\in[S_1,S_2], k\in K_u, w_i\in[-W,W]$. Moreover,
 for any $t,t'\in [S_1,S_2]$ uniformly with respect to $k\in K_u$, any $w_i\in \R$,
\BQN\label{XVar}
\Var\LT((1+b_iu^{-2}f_{i}(u^{\lambda_{i}}t))\mathcal{X}^w_{u,i}(t,k)-(1+b_iu^{-2}
f_{i}(u^{\lambda_{i}}t'))\mathcal{X}^w_{u,i}(t',k)\RT)\rw 2c_i^2a_i|t-t'|^\alpha\II_{\{\alpha_i=\alpha\}},
\EQN
as $u\rw\IF$.
Combination of (\ref{XM12}) and (\ref{XVar}) shows that the finite-dimensional distributions of
$$\{(\vk{1}+\vk{b}u^{-2}\vk{f}(u^{\vk{\lambda}}t))\vk{\mathcal{X}}^{\vk{w}}_u(t,k), t\in [S_1, S_2]\}$$ weakly converge to the finite-dimensional
distributions of  $\{\vk{\zeta}(t), t\in [S_1, S_2]\}$. Moreover,
by (\ref{rp}) we have that there exists a constant $C>0$ such that for all $t, t'\in[S_1,S_2]$ and  all large  $u$
\BQN\label{HC}
&&\sup_{k\in K_u}\Var\LT((1+b_iu^{-2}f_{i}(u^{\lambda_{i}}t))
\mathcal{X}^w_{u,i}(t,k)-(1+b_iu^{-2}f_{i}(u^{\lambda_{i}}t'))
\mathcal{X}^w_{u,i}(t',k)\RT)\nonumber\\
&&\quad \leq m_{u,i}^2(k)
Var\LT(\xi_{i}(u^{-2/\alpha}t)-\xi_{i}(u^{-2/\alpha}t')\RT)\leq C|t-t'|^\alpha,
\EQN
which combined with \eqref{XM12} implies that the family of distributions
$$\pk{(\vk{1}+\vk{b}u^{-2}\vk{f}(u^{\vk{\lambda}}t))\vk{\mathcal{X}}^{\vk{w}}_u(t,k)\in(\cdot)}$$
is uniformly tight with respect to $k\in K_u$ and $\vk{w}$ in a compact set of $\R^n$. Consequently,
$$\{(\vk{1}+\vk{b}u^{-2}\vk{f}(u^{\vk{\lambda}}t))\vk{\mathcal{X}}^{\vk{w}}_u(t,k), t\in [S_1,S_2]\} \quad \text{weakly converges to} \quad \{\vk{\zeta}(t), t\in [S_1,S_2]\}.$$
Since
$$\lim_{u\rw\IF}\max_{1\leq i\leq n}\sup_{k\in K_u}\sup_{t\in [S_1,S_2]}\left|(1+b_iu^{-2}f_i(u^{\lambda_i}t))-1\right|=0,$$
we conclude that
$$\{\vk{\mathcal{X}}^{\vk{w}}_u(t,k), t\in [S_1,S_2]\} \quad \text{weakly converges to} \quad \{\vk{\zeta}(t), t\in [S_1,S_2]\}.$$
\underline{Upper bound for $I_1(u)$}. We first show that
\BQNY
c_u(\vk{w}):&=&\sup_{k\in K_u}\LT|\pk{\exists_{t\in
[S_1,S_2]}\vk{\mathcal{X}}^{\vk{w}}_u(t,k)>\vk{w}}
-\pk{\exists_{t\in[S_1,S_2]}\vk{\zeta}(t)>\vk{w}}\RT|\\
&=&
\sup_{k\in K_u}\LT|\pk{\sup_{t\in[S_1,S_2]}\min_{1\leq i\leq n}(\mathcal{X}^{w}_{u,i}(t,k)-w_i)>0}
-\pk{\sup_{t\in[S_1,S_2]}\min_{1\leq i\leq n}(\zeta_i(t)-w_i)>0}\RT|\rw 0,
\EQNY
for almost all $\vk{w}\in \mathbb{R}^n$.
Let $$\mathbb{A}:=\left\{\vk{v}: \pk{\sup_{t\in[S_1,S_2]}\min_{1\leq i\leq n}(\zeta_i(t)-v_i)>0}  \text{is continuous at }  \vk{v}\right\}.$$
Note that if $\vk{w}\in \mathbb{A}$, then
 $$\pk{\sup_{t\in[S_1,S_2]}\min_{1\leq i\leq n}(\zeta_i(t)-w_i)>x}$$
 is continuous with respect to $x$ at $x=0$.
Hence by the continuity of functional $\sup\min$, we have that
\BQNY
c_u(\vk{w})\rw 0,
\EQNY
for $\vk{w}\in \mathbb{A}$ and  $mes(\mathbb{A}^c)=0$.
Thus in light of dominated convergence theorem, we have
 \BQNY
 I_1(u)\leq e^{nW}\int_{\vk{w}\in \widetilde{W}^n \cap \mathbb{A} }c_u(\vk{w})d\vk{w}+W^n e^{nW}\sup_{\vk{w}\in \widetilde{W}^n}\left|1-e^{-\sum_{i=1}^n\frac{w_i^2}{2m_{u,i}^2(k)}}\right|\rw 0, \quad u\rw\IF.
 \EQNY
\underline{Upper bound for $I_2(u)$}.
Using (\ref{XM11}) and (\ref{XVar}), for some $\delta\in(0,1/2)$, $|w_i|>W$ with $W$ sufficiently large and all $u$ large we have
\BQNY
\sup_{k\in K_u,t\in[S_1,S_2]}\EF{(1+b_iu^{-2}f_{i}(u^{\lambda_{i}}t))\mathcal{X}^w_{u,i}(t,k)}
\leq\mathbb{C}_1+\delta|w_i|
\EQNY
and
\BQNY
\sup_{k\in K_u,t\in[S_1,S_2]}\Var\LT((1+b_iu^{-2}f_{i}(u^{\lambda_{i}}t))\mathcal{X}^w_{u,i}(t,k)\RT)
\leq \mathbb{C}_2.
\EQNY
\COM{
\BQNY
\sup_{k\in K_u,t\in[S_1,S_2]}\Var\LT(\sum_{i=1}^n\mathcal{X}^w_{u,i}(t,k)\RT)\leq \mathbb{C}_4.
\EQNY}
Moreover, by the mutual independence of $\vk{\mathcal{X}}^w_{u,i}(t,k), 1\leq i\leq n$
\BQNY
\pk{\exists_{t\in
[S_1,S_2]}\vk{\mathcal{X}}^{\vk{w}}_u(t,k)>\vk{w} }&=&\pk{\sup_{t\in
[S_1,S_2]}\min_{1\leq i\leq n}\LT({\mathcal{X}}^w_{u,i}(t,k)-w_i\RT)>0}\\
&\leq& \pk{\min_{1\leq i\leq n}\LT(\sup_{t\in
[S_1,S_2]}{\mathcal{X}}^w_{u,i}(t,k)-w_i\RT)>0}\\
&=&\prod_{i=1}^n \pk{\sup_{t\in
[S_1,S_2]}{\mathcal{X}}^w_{u,i}(t,k)>w_i}.
\EQNY
Consequently, it follows that
\BQNY
\sup_{k\in K_u}\int_{\widetilde{W}_j^n}e^{\sum_{i=1}^nw_i}\pk{\exists_{t\in
[S_1,S_2]}\vk{\mathcal{X}}^{\vk{w}}_u(t,k)>\vk{w} }d\vk{w}\leq J_1\times J_2,
\EQNY
where by (\ref{HC}) and Theorem 8.1 of \cite{Pit96}
\BQNY
J_1&=&\sup_{k\in K_u}\int_{|w_j|>W}e^{w_j}\pk{\sup_{t\in
[S_1,S_2]}\mathcal{X}^w_{u,j}(t,k)>w_j}d w_j\\
&\leq&\sup_{k\in K_u}\int_{|w_j|>W}e^{w_j}\mathbb{P}\left(\sup_{t\in
[S_1,S_2]}\LT((1+b_iu^{-2}f_{i}(u^{\lambda_{i}}t))\mathcal{X}^w_{u,j}(t,k)
-\EF{(1+b_iu^{-2}f_{i}(u^{\lambda_{i}}t))\mathcal{X}^w_{u,j}(t,k)}\RT)\right.\\
&&\quad \left.>(1-\delta)|w_j|-\mathbb{C}_1\right)d w_j\\
&\leq& e^{-W}+\int_W^{\IF}e^{w_j}\mathbb{C}_3{w_j}^{2/\alpha}
\Psi\LT(\frac{(1-\delta)w_j-\mathbb{C}_1}{\mathbb{C}_2}\RT)d w_j\\
&=:&A_1(W)\rw 0,\ W\rw\IF,
\EQNY
and
{\small
\BQNY
J_2&=&\sup_{k\in K_u}\underset{i\neq j}{\prod_{i=1}^n}\LT(\int_{\R}e^{w_i}
\pk{\sup_{t\in
[S_1,S_2]}{\mathcal{X}}^w_{u,i}(t,k)>w_i}d w_i\RT)\\
&\leq&\sup_{k\in K_u}\underset{i\neq j}{\prod_{i=1}^n}\LT(e^{W_1}+\right.\\
&&\left.+\int_{W_1}^\IF e^{w_i} \mathbb{P}\left(\sup_{t\in
[S_1,S_2]}\LT((1+b_iu^{-2}f_{i}(u^{\lambda_{i}}t)){\mathcal{X}}^w_{u,i}(t,k)
-\EF{(1+b_iu^{-2}f_{i}(u^{\lambda_{i}}t))\mathcal{X}^w_{u,i}(t,k)}\RT)
>(1-\delta)w_i-\mathbb{C}_1\right)d w_i\RT)\\
&\leq&\underset{i\neq j}{\prod_{i=1}^n}\LT(e^{W_1}+\int_{W_1}^\IF e^{w_i}\mathbb{C}_4{w_i}^{2/\alpha}
\Psi\LT(\frac{(1-\delta)w_i-\mathbb{C}_1}{\mathbb{C}_2}\RT)d w_i\RT)
\leq \mathbb{C}_5,
\EQNY
}
with $W_1$  some positive constant.
Thus we have
\BQNY
I_2(u)\leq n\mathbb{C}_5 A_1(W)\rw 0,\ W\rw \IF.
\EQNY
\underline{Upper bound for $I_3(u)$}.
Borell-TIS inequality (see, e.g., \cite{AdlerTaylor}) implies that  $$I_3(u)\rw 0,\ u, W\rw\IF.$$
Hence (\ref{IR}) follows.\\
ii)  Suppose that $\lambda>2/\alpha$.
Observe that
\BQNY
&&\frac{\pk{\exists_{t\in
[u^{-\lambda}S_1,u^{-\lambda}S_2]}\vk{Z}_u(t)>\vk{m}_u(k)}}{\prod_{i=1}^n\Psi(m_{u,i}(k))}\\
&&=\LT(\prod_{i=1}^n\frac{e^{-\frac{\left(m_{u,i}(k)\right)^2}{2}}}{\sqrt{2\pi}m_{u,i}(k)
\Psi(m_{u,i}(k))}\RT)\int_{{\R}^n}
e^{\sum_{i=1}^n\left(w_i-\frac{w_i^2}{2\left(m_{u,i}(k)\right)^2}\right)}
\pk{\exists_{t\in
[S_1,S_2]}\vk{\mathcal{X}}^{\vk{w}}_u(t,k)>\vk{w} }d\vk{w},
\EQNY
where $\vk{\mathcal{X}}^{\vk{w}}_u(t,k)=(\mathcal{X}^w_{u,1}(t,k),\ldots,\mathcal{X}^w_{u,n}(t,k))$ with
\BQNY
\mathcal{X}^{w}_{u,i}(t,k)=m_{u,i}(k)(Z_{u,i}(u^{-\lambda}t)-m_{u,i}(k))+w_i
\Bigl\lvert\xi_i(0)=
m_{u,i}(k)-\frac{w_i}{m_{u,i}(k)}.
\EQNY
The rest of derivations for this case is the same as given in the proof for case $\lambda\leq 2/\alpha$, with
exception that
\BQNY
\EF{(1+b_iu^{-2}f_{i}(u^{\lambda_{i}}t))\mathcal{X}^w_{u,i}(t,k)}\rw
 -c_i^2\widetilde{f}_{i}(t),\ u\rw\IF,
\EQNY
and
\BQNY
\Var\LT(\mathcal{X}^w_{u,i}(t,k)-\mathcal{X}^w_{u,i}(t',k)\RT)\rw 0,\ u\rw\IF.
\EQNY
Hence we omit the rest of the proof.
\QED
\BEL\label{pan3} Let $\vk{X}(t), (t)\in\R$ be a centered vector-valued stationary Gaussian process
 with independent coordinates $X_i$'s.
Suppose that for each $i=1,...,n$,  $X_i(t)$ has continuous sample paths, unit variance and correlation function $r_i(\cdot),\ 1\leq i\leq n,$ satisfying
\BQN\label{rp1}
0<1-2a_i|t|^{\alpha_i}\leq r_i(t)\leq1-\frac{a_i}{2}|t|^{\alpha_i},\ a_i>0, \ \alpha_i\in(0,2],
\EQN
for all $t\in[0,\vn]$  with $0<\vn<1$ small enough. Let $K_u$ be an index set. Then
we have for any \cLa{$\vk{m}_{u}(k)$, $\vk{w}_{u}(l)$} such that
\BQN\label{Mku2}
\lim_{u\to\IF}\sup_{k\in K_u}\LT|\frac{1}{u}\vk{m}_{u}(k)-\vk{c}\RT|=0,\ \lim_{u\to\IF}\sup_{l\in K_u}\LT|\frac{1}{u}\vk{w}_{u}(l)-\vk{c}\RT|=0,
\EQN
and any $T(k,l)>S>1$ satisfying $\lim_{u\rw\IF}\sup_{k,l\in K_u}\frac{T(k,l)}{u^{2/\alpha}}=0$, that
\BQNY
&&\pk{\exists_{t\in[0,S]u^{-2/\alpha}}\vk{X}(t)>\vk{m}_u(k),\exists_{t\in[T(k,l), T(k,l)+S]u^{-2/\alpha}}\vk{X}(t)>\vk{w}_u(l)}\\
&&\quad\quad\quad\quad \ \ \ \ \ \ \  \ \ \ \ \ \ \leq FS^{2n}\exp(-G(T(k,l)-S)^{\alpha})\prod_{i=1}^n\Psi\LT(\frac{m_{u,i}(k)+w_{u,i}(l)}{2}\RT)
\EQNY
holds uniformly for any $k, l \in K_u$ and all $u$ large
where $\alpha=\min_{1\leq i\leq n}(\alpha_i) $ and $F,G$ are two positive constants.
\EEL
\prooflem{pan3}
By the independence of $X_i$'s, we have that
\BQNY
&&\pk{\exists_{t\in[0,S]u^{-2/\alpha}}\vk{X}(t)>\vk{m}_u(k),\exists_{t\in[T(k,l), T(k,l)+S]u^{-2/\alpha}}\vk{X}(t)>\vk{w}_u(l)}\\
&&\quad\quad\leq\pk{\bigcap_{i=1}^n\LT\{\sup_{t\in[0,S]u^{-2/\alpha}}X_i(t)>m_{u,i}(k)\RT\},
\bigcap_{i=1}^n\LT\{\sup_{t\in[T(k,l), T(k,l)+S]u^{-2/\alpha}}X_i(t)>w_{u,i}(k)\RT\}}\\
&&\quad\quad\bl{\leq}\prod_{i=1}^n\pk{\sup_{t\in[0,S]u^{-2/\alpha}}X_i(t)>m_{u,i}(k),
\sup_{t\in[T(k,l), T(k,l)+S]u^{-2/\alpha}}X_i(t)>w_{u,i}(k)}.
\EQNY
Application of Lemma 6.3 in \cite{Pit96} (or Theorem 3.1 in \cite{Uniform2016}) for each term in the above product establishes the claim.
\QED

\prooftheo{PreThm1}
Let
\BQNY
\pi(u):=\pk{\exists_{t\in\E} \vk{X}_{u}(t)>\vk{m}_u}=\pk{ \exists_{t\in\E} \frac{\vk{X}_u(t)}{\vk{\sigma}_u(t)}\frac{\vk{\sigma}_u(t)}{\vk{\sigma}_u(0)}>
\frac{\vk{m}_u}{\vk{\sigma}_u(0)}}.
\EQNY
In view of {\bf A2-A3} and by Gordon inequality (see, e.g., Lemma 5.1 in \cite{Tabis}), we have that for $\epsilon \in (0,1)$ and $u$ sufficiently large
 \BQN\label{main}
\pk{\exists_{t\in\E} \vk{Z}_{u,-\vn}(t)>\frac{\vk{m}_u}{\vk{\sigma}_u(0)}}\leq\pi(u)\leq\pk{\exists_{t\in\E} \vk{Z}_{u,+\vn}(t)>\frac{\vk{m}_u}{\vk{\sigma}_u(0)}}.
\EQN
where
$$\vk{Z}_{u,\pm\vn}(t)=\frac{\vk{Y}_{\pm\vn}(t)}{\vk{w}_{u,\mp\vn}(t)},\ t\inr,$$
with  $\vk{Y}_{\pm\vn}(t),t\inr$ being  homogeneous vector-valued Gaussian processes with independent coordinates $Y_{i,\pm\vn}(t),t\inr$ having continuous trajectories, unit variance and correlation function satisfying
$$r_{i,\pm\vn}(t)=e^{-(1\pm\vn)a_i|t|^{\alpha_i}},$$
and
\cLa{$\vk{w}_{u,\pm\vn}(t)=(w_{u,1,\pm\vn}(t),\ldots,w_{u,n,\pm\vn}(t))$ with
\BQNY
w_{u,i,\pm\vn}(t)=1+u^{-2}\LT(f _{i}(u^{\lambda_{i}}t)\pm\vn\abs{f_{i}(u^{\lambda_{i}}t)}\pm\vn\RT), \quad \epsilon\in (0,1).
\EQNY}
Next, we use the double-sum method to derive an upper and a lower bound of (\ref{main}) and then show that they are asymptotically tight.
We distinguish three scenarios: $\lambda<2/\alpha$, $\lambda=2/\alpha$ and $\lambda<2/\alpha$.\\
{\it $\diamond$ \underline{Case $\lambda<2/\alpha$}.}  For any $S>0$, let
\BQN\label{notation1}
I_k(u)&=&[ku^{-2/\alpha}S,(k+1)u^{-2/\alpha}S], \ \ k\in\mathbb{Z},\ N_1(u)=\LT\lfloor \frac{x_1(u)}{Su^{-2/\alpha}}\RT\rfloor-\mathbb{I}_{\{x_1\leq0\}},\nonumber\\ N_2(u)&=&\LT\lfloor \frac{x_2(u)}{Su^{-2/\alpha}}\RT\rfloor+\mathbb{I}_{\{x_2\leq0\}}, \quad \cLa{\vk{v}_{u,\pm\vn}(k)=(v_{u,1,\pm\vn}(k),\ldots,v_{u,n,\pm\vn}(k)),}
\EQN
 with
\BQNY
v_{u,i,+\vn}(k)=\frac{m_{u,i}}{\sigma_{u,i}(0)}\sup_{s\in I_k(u)}w_{u,i,+\vn}(s),\ \
v_{u,i,-\vn}(k)=\frac{m_{u,i}}{\sigma_{u,i}(0)}\inf_{s\in I_k(u)}w_{u,i,-\vn}(s).
\EQNY
For $u$ large enough, in view of (\ref{main}) we have
\BQNY
&&\pi(u)\leq\pk{\exists_{t\in\E} \vk{Z}_{u,+\vn}(t)>\frac{\vk{m}_u}{\vk{\sigma}_u(0)}}\leq\sum_{k=N_1(u)}^{N_2(u)}\pk{ \exists_{t\in I_{k}(u)} \vk{Z}_{u,+\vn}(t)>\frac{\vk{m}_u}{\vk{\sigma}_u(0)}},\label{proof31}\\
&&\pi(u)\geq\pk{\exists_{t\in\E} \vk{Z}_{u,-\vn}(t)>\frac{\vk{m}_u}{\vk{\sigma}_u(0)}}\geq\sum_{k=N_1(u)+1}^{N_2(u)-1}
\pk{ \exists_{t\in I_{k}(u)} \vk{Z}_{u,-\vn}(t)>\frac{\vk{m}_u}{\vk{\sigma}_u(0)}}-\sum_{i=1}^2\Lambda_i(u),\label{proof32}
\EQNY
where
$$\Lambda_1(u)=\sum_{k=N_1(u)}^{N_2(u)}\pk{ \exists_{t\in I_{k}(u)} \vk{Z}_{u,-\vn}(t)>\frac{\vk{m}_u}{\vk{\sigma}_u(0)}, \exists_{t\in I_{k+1}(u)} \vk{Z}_{u,-\vn}(t)>\frac{\vk{m}_u}{\vk{\sigma}_u(0)}},$$
and
$$\Lambda_2(u)=\sum_{N_1(u)\leq k,l\leq N_2(u), l\geq k+2}\pk{ \exists_{t\in I_{k}(u)} \vk{Z}_{u,-\vn}(t)>\frac{\vk{m}_u}{\vk{\sigma}_u(0)}, \exists_{t\in I_{l}(u)} \vk{Z}_{u,-\vn}(t)>\frac{\vk{m}_u}{\vk{\sigma}_u(0)}}.$$
{\it \underline{Asymptotics of $\pi(u)$}}.
By stationarity of $\vk{Y}_{+\epsilon}$ and \nelem{pan2}, we have that
\BQNY
\pi(u)
&\leq& \sum_{k=N_1(u)}^{N_2(u)}\pk{\exists_{t\in I_{k}(u)} \vk{Y}_{+\vn}(t)>\vk{v}_{u,-\vn}(k)}\nonumber\\
&\leq& \sum_{k=N_1(u)}^{N_2(u)}\pk{\exists_{t\in I_{0}(u)} \vk{Y}_{+\vn}(t)>\vk{v}_{u,-\vn}(k)}\nonumber\\
&\sim& \mathcal{H}_{\alpha,(1+\vn)\frac{\vk{a}}{{\vk{\sigma}}^2}
\mathbf{I}_{\{\vk{\alpha}=\alpha\vk{1}\}}}[0,S]
\sum_{k=N_1(u)}^{N_2(u)}\prod_{i=1}^n\Psi(v_{u,i,-\vn}(k)),\ u\rw\IF.
\EQNY
Furthermore,
\BQN\label{regular}
&&\sum_{k=N_1(u)}^{N_2(u)}\prod_{i=1}^n\Psi(v_{u,i,-\vn}(k))\nonumber\\
&&\sim\sum_{k=N_1(u)}^{N_2(u)}\prod_{i=1}^n\LT(\frac{1}{\sqrt{2\pi}v_{u,i,-\vn}(k)}
\exp\LT(-\frac{v^2_{u,i,-\vn}(k)}{2}\RT)\RT)\nonumber\\
&&\sim\LT(\prod_{i=1}^n\Psi\LT(\frac{m_{u,i}}{\sigma_{u,i}(0)}\RT)\RT)
\sum_{k=N_1(u)}^{N_2(u)}\exp\LT(-\sum_{i=1}^n\frac{m_{u,i}^2u^{-2}\inf_{s\in I_k(u)}\LT(f_{i}(u^{\lambda_{i}}s)-\vn\abs{f_{i}(u^{\lambda_{i}}s)}-\vn\RT)}
{\sigma^2_{u,i}(0)}\RT)\nonumber\\
&&\sim\LT(\prod_{i=1}^n\Psi\LT(\frac{m_{u,i}}{\sigma_{u,i}(0)}\RT)\RT)
\sum_{k=N_1(u)}^{N_2(u)}\exp\LT(-\sum_{i=1}^n\frac{m_{u,i}^2u^{-2}\inf_{s\in [k,k+1]}\LT(f_{i}(u^{\lambda_{i}-\frac{2}{\alpha}}S s)-\vn\abs{f_{i}(u^{\lambda_{i}-\frac{2}{\alpha}} S s)}-\vn\RT)}{\sigma^2_{u,i}(0)}\RT)\nonumber\\
&&\leq\LT(\prod_{i=1}^n\Psi\LT(\frac{m_{u,i}}{\sigma_{u,i}(0)}\RT)\RT)S^{-1}u^{2/\alpha-\lambda}
\int_{x_1}^{x_2}\exp\LT(-\sum_{i=1}^n\frac{\widetilde{f}_i^{\vn}(t)}{\sigma_i^2}\RT)dt,
\EQN
where $\widetilde{f}_i^{\vn}(t)=\widetilde{f}_{i}\LT(t\RT)-\vn \abs{\widetilde{f}_{i}\LT(t\RT)}-\vn$. \bl{ In order to prove (\ref{regular}), we note that for $-\IF<x_1<x_2<\IF$,
\BQNY
&&\sum_{k=N_1(u)}^{N_2(u)}\exp\LT(-\sum_{i=1}^n\frac{m_{u,i}^2u^{-2}\inf_{s\in [k,k+1]}\LT(f_{i}(u^{\lambda_{i}-\frac{2}{\alpha}}S s)-\vn\abs{f_{i}(u^{\lambda_{i}-\frac{2}{\alpha}} S s)}-\vn\RT)}{\sigma^2_{u,i}(0)}\RT)\\
&&\sim S^{-1}u^{2/\alpha-\lambda}
\int_{x_1}^{x_2}\exp\LT(-\sum_{i=1}^n\frac{\widetilde{f}_i^{\vn}(t)}{\sigma_i^2}\RT)dt,\quad u\rw\IF,
\EQNY
which implies that (\ref{regular}) holds for $-\IF<x_1<x_2<\IF$.  Next we assume that $-\IF<x_1<x_2=\IF$. Let $y$ be a positive constant satisfying $x_1<y<\IF$ and $N(u,y)=\left[\frac{yu^{2/\alpha-\lambda}}{S}\right]$. Then it follows that
\BQN\label{first}
&&\sum_{k=N_1(u)}^{N(u,y)}\exp\LT(-\sum_{i=1}^n\frac{m_{u,i}^2u^{-2}\inf_{s\in [k,k+1]}\LT(f_{i}(u^{\lambda_{i}-\frac{2}{\alpha}}S s)-\vn\abs{f_{i}(u^{\lambda_{i}-\frac{2}{\alpha}} S s)}-\vn\RT)}{\sigma^2_{u,i}(0)}\RT)\nonumber\\
&&\sim S^{-1}u^{2/\alpha-\lambda}
\int_{x_1}^{y}\exp\LT(-\sum_{i=1}^n\frac{\widetilde{f}_i^{\vn}(t)}{\sigma_i^2}\RT)dt, \quad u\rw\IF.
\EQN
 By Potter's Theorem (Theorem 1.5.6 in \cite{BI1989}) and the fact that for $j\in\{1\leq i\leq n: \lambda_i=\lambda\}$, $f_j(t)$ is regularly varying at $\IF$ with positive index, we have that for any $\eta>0$ and sufficiently large $y$ and $u$
\BQNY \left|\frac{\sigma^2_{j}}{\sigma^2_{u,j}(0)}\frac{m_{u,j}^2u^{-2}\inf_{s\in [k,k+1]}\LT(f_{j}(u^{\lambda-\frac{2}{\alpha}}S s)-\vn\abs{f_{j}(u^{\lambda-\frac{2}{\alpha}} S s)}-\vn\RT)}{\widetilde{f}_j^{\vn}(u^{\lambda-\frac{2}{\alpha}}S k)}-1\right|<\eta
\EQNY
holds  for all $k> N(u,y)$.
Then we have that for $k> N(u,y)$
\BQNY
\left|\sum_{\lambda_i=\lambda}\frac{m_{u,i}^2u^{-2}\inf_{s\in [k,k+1]}\LT(f_{i}(u^{\lambda_{i}-\frac{2}{\alpha}}S s)-\vn\abs{f_{i}(u^{\lambda_{i}-\frac{2}{\alpha}} S s)}-\vn\RT)}{\sigma^2_{u,i}(0)}-\sum_{\lambda_i=\lambda}
\frac{\widetilde{f}_i^{\vn}(u^{\lambda-\frac{2}{\alpha}}S k)}{\sigma_i^2}\right|\leq \eta \sum_{\lambda_i=\lambda}\frac{|\widetilde{f}_i^{\vn}(u^{\lambda-\frac{2}{\alpha}}S k)|}{\sigma_i^2}
\EQNY
Using (\ref{xx}), it follows that
\BQNY
\lim_{u\rw\IF}\sup_{N_1(u)\leq k\leq N_2(u)}\left|\sum_{\lambda_i<\lambda}\frac{m_{u,i}^2u^{-2}\inf_{s\in [k,k+1]}\LT(f_{i}(u^{\lambda_{i}-\frac{2}{\alpha}}S s)-\vn\abs{f_{i}(u^{\lambda_{i}-\frac{2}{\alpha}} S s)}-\vn\RT)}{\sigma^2_{u,i}(0)}-\sum_{\lambda_i<\lambda}
\frac{\widetilde{f}_i^{\vn}(u^{\lambda_i-\frac{2}{\alpha}}S k)}{\sigma_i^2}\right|=0.
\EQNY
Hence,   for sufficiently large $y$ and $u$  we have that
\BQNY
\sum_{i=1}^n\frac{m_{u,i}^2u^{-2}\inf_{s\in [k,k+1]}\LT(f_{i}(u^{\lambda_{i}-\frac{2}{\alpha}}S s)-\vn\abs{f_{i}(u^{\lambda_{i}-\frac{2}{\alpha}} S s)}-\vn\RT)}{\sigma^2_{u,i}(0)}\geq \sum_{i=1}^n\frac{\widetilde{f}_i^{\vn}(u^{\lambda-\frac{2}{\alpha}}S k)}{\sigma_i^2}-\eta\sum_{i=1}^n\frac{|\widetilde{f}_i^{\vn}
(u^{\lambda-\frac{2}{\alpha}}S k)|}{\sigma_i^2}
\EQNY
holds for $k> N(u,y)$.
Combining the above with \eqref{FF} implies that
\BQNY
&&\sum_{k=N(u,y)+1}^{N_2(u)}\exp\LT(-\sum_{i=1}^n\frac{m_{u,i}^2u^{-2}\inf_{s\in [k,k+1]}\LT(f_{i}(u^{\lambda_{i}-\frac{2}{\alpha}}S s)-\vn\abs{f_{i}(u^{\lambda_{i}-\frac{2}{\alpha}} S s)}-\vn\RT)}{\sigma^2_{u,i}(0)}\RT)\\
&&\leq \sum_{k=N(u,y)+1}^{N_2(u)}\exp\LT(-\sum_{i=1}^n\frac{\widetilde{f}_i^{\vn}(u^{\lambda-\frac{2}{\alpha}}S k)}{\sigma_i^2}+\eta\sum_{i=1}^n\frac{|\widetilde{f}_i^{\vn}(u^{\lambda-\frac{2}{\alpha}}S k)|}{\sigma_i^2}\RT)\\
&&\leq u^{\frac{2}{\alpha}-\lambda}S^{-1}\int_{y}^{\IF}\exp\LT(-\sum_{i=1}^n\frac{\widetilde{f}_i^{\vn}(t)}{\sigma_i^2}+
\eta\sum_{i=1}^n\frac{|\widetilde{f}_i^{\vn}(t)|}{\sigma_i^2}\RT)dt,
\EQNY
which together with (\ref{first}) and  the arbitrariness of $\eta>0$ confirms that (\ref{regular}) holds. For other cases of $x_1$ and $x_2$, we can similarly show that (\ref{regular}) is satisfied.}
 By (\ref{xx}) and \eqref{FF}, we have that $$\int_{x_1}^{x_2}\exp\LT(-\sum_{i=1}^n\frac{\widetilde{f}_i^{\vn}(t)}{\sigma_i^2}\RT)dt<\IF.$$
Consequently,
\BQN\label{proof4}
\pi(u)\leq\mathcal{H}_{\alpha,\frac{\vk{a}}{{\vk{\sigma}}^2}\mathbf{I}_{\{\vk{\alpha}=\alpha\vk{1}\}}}
u^{2/\alpha-\lambda}\int_{x_1}^{x_2}\exp\LT(-\sum_{i=1}^n\frac{\widetilde{f}_i(t)}{\sigma_i^2}\RT)dt
\LT(\prod_{i=1}^n\Psi\LT(\frac{m_{u,i}}{\sigma_{u,i}(0)}\RT)\RT),
\EQN
as $u\rw\IF,\ S\rw\IF,\ \vn\rw 0$.
Analogously, we have
\BQN\label{proof5}
&&\sum_{k=N_1(u)+1}^{N_2(u)-1}\pk{ \exists_{t\in I_{k}(u)} \vk{Z}_{u,-\vn}(t)>\frac{\vk{m}_u}{\vk{\sigma}_u(0)}}\nonumber\\
&&\quad\quad\geq\mathcal{H}_{\alpha,\frac{\vk{a}}{{\vk{\sigma}}^2}
\mathbf{I}_{\{\vk{\alpha}=\alpha\vk{1}\}}}
u^{2/\alpha-\lambda}\int_{x_1}^{x_2}\exp\LT(-\sum_{i=1}^n\frac{\widetilde{f}_i(t)}{\sigma_i^2}\RT)dt
\LT(\prod_{i=1}^n\Psi\LT(\frac{m_{u,i}}{\sigma_{u,i}(0)}\RT)\RT),
\EQN
as $u\rw\IF,\ S\rw\IF,\ \vn\rw 0$.\\
{\it \underline{Upper bound for $\Lambda_1(u)$}}.
It follows that
\BQN\label{proof6}
\Lambda_1(u)&=&\sum_{k=N_1(u)}^{N_2(u)}\left(\pk{\exists_{t\in I_{k}(u)} \vk{Z}_{u,-\vn}(t)>\frac{\vk{m}_u}{\vk{\sigma}_u(0)}}+\pk{\exists_{t\in I_{k+1}(u)} \vk{Z}_{u,-\vn}(t)>\frac{\vk{m}_u}{\vk{\sigma}_u(0)}}\RT.\nonumber\\
& &\LT.-\pk{\exists_{t\in I_{k}(u)\cup I_{k+1}(u)} \vk{Z}_{u,-\vn}(t)>\frac{\vk{m}_u}{\vk{\sigma}_u(0)}}\right)\nonumber\\
&\leq&\sum_{k=N_1(u)}^{N_2(u)}\left(\pk{\exists_{t\in I_{k}(u)} \vk{Y}_{-\vn}(t)>\vk{\widehat{v}}_{u,+\vn}(k)}+\pk{\exists_{t\in I_{k+1}(u)} \vk{Y}_{-\vn}(t)>\vk{\widehat{v}}_{u,+\vn}(k)}\RT.\nonumber\\
& &\LT.-\pk{\exists_{t\in I_{k}(u)\cup I_{k+1}(u)} \vk{Y}_{-\vn}(t)>\vk{\widetilde{v}}_{u,+\vn}(k)}\right)\nonumber\\
&\sim&\left(2\mathcal{H}_{\alpha,(1-\vn)\frac{\vk{a}}{{\vk{\sigma}}^2}
\mathbf{I}_{\{\vk{\alpha}=\alpha\vk{1}\}}}[0,S]
-\mathcal{H}_{\alpha,(1-\vn)\frac{\vk{a}}{{\vk{\sigma}}^2}\mathbf{I}_{\{\vk{\alpha}
=\alpha\vk{1}\}}}[0,2S]\right)
\sum_{k=N_1(u)}^{N_2(u)}\LT(\prod_{i=1}^n\Psi(v_{u,i,+\vn}(k))\RT)\nonumber\\
&=&o\left(u^{2/\alpha-\lambda}
\prod_{i=1}^n\Psi\LT(\frac{m_{u,i}}{\sigma_{u,i}(0)}\RT)\right), \ u\rw\IF, S\rw \IF, \vn\rw 0,
\EQN
where $$\widehat{v}_{u,i,+\vn}(k)=\min \LT(\frac{m_{u,i}}{\sigma_{u,i}(0)}\inf_{s\in I_k(u)}w_{u,i,+\vn}(s),\frac{m_{u,i}}{\sigma_{u,i}(0)}\inf_{s\in I_{k+1}(u)}w_{u,i,+\vn}(s)\RT)$$
and
$$\widetilde{v}_{u,i,+\vn}(k)=\max\LT(v_{u,i,+\vn}(k),
v_{u,i,+\vn}(k+1)\RT).$$
{\it\underline{Upper bound for $\Lambda_2(u)$}}.
In light of \nelem{pan3}, we have that
\BQN\label{proof7}
\Lambda_2(u)&=& \sum_{N_1(u)\leq k,l\leq N_2(u), l\geq k+2}\pk{ \exists_{t\in I_{k}(u)} \vk{Z}_{u,-\vn}(t)>\frac{\vk{m}_u}{\vk{\sigma}_u(0)}, \exists_{t\in I_{l}(u)} \vk{Z}_{u,-\vn}(t)>\frac{\vk{m}_u}{\vk{\sigma}_u(0)}}\nonumber\\
&\leq& \sum_{N_1(u)\leq k,l\leq N_2(u), l\geq k+2}\pk{\exists_{t\in I_{k}(u)} \vk{Y}_{-\vn}(t)>\vk{\overline{v}}_{u,+\vn}(k),\exists_{t\in I_{l}(u)} \vk{Y}_{-\vn}(t)>\vk{\overline{v}}_{u,+\vn}(l)}\nonumber\\
&\leq& \sum_{N_1(u)\leq k,l\leq N_2(u), l\geq k+2}\pk{\exists_{t\in I_{0}(u)} \vk{Y}_{-\vn}(t)>\vk{\overline{v}}_{u,+\vn}(k),\exists_{t\in I_{l-k}(u)} \vk{Y}_{-\vn}(t)>\vk{\overline{v}}_{u,+\vn}(l)}\nonumber\\
&\leq&\sum_{N_1(u)\leq k,l\leq N_2(u), l\geq k+2}\mathbb{C}_1S^{2n}\exp(-\mathbb{C}_2((l-k-1)S)^\alpha)
\prod_{i=1}^n\Psi\LT(\frac{\overline{v}_{u,i,-\vn}(k)+\overline{v}_{u,i,-\vn}(l)}{2}\RT) \nonumber\\
&\leq&2\sum_{l=1}^\IF \mathbb{C}_1S^{2n}\exp(-\mathbb{C}_2(l S)^\alpha)
\sum_{k=N_1(u)}^{N_2(u)}\prod_{i=1}^n\Psi\LT(\overline{v}_{u,i,-\vn}(k)\RT) \nonumber\\
&\leq&S^{2n}\exp(-\mathbb{C}_3S^\alpha)
u^{2/\alpha-\lambda}
\prod_{i=1}^n\Psi\LT(\frac{m_{u,i}}{\sigma_{u,i}(0)}\RT) \nonumber\\
&=&o\left(u^{2/\alpha-\lambda}
\prod_{i=1}^n\Psi\LT(\frac{m_{u,i}}{\sigma_{u,i}(0)}\RT)\right),  \ u\rw\IF, S\rw\IF,
\EQN
where
$$\overline{v}_{u,i,+\vn}(k)=\frac{m_{u,i}}{\sigma_{u,i}(0)}\inf_{s\in I_k(u)}w_{u,i,+\vn}(s).$$
Combination of  (\ref{proof31})-(\ref{proof7}) leads to
$$\pi(u)\sim\mathcal{H}_{\alpha,\frac{\vk{a}}{{\vk{\sigma}}^2}\mathbf{I}_{\{\vk{\alpha}=\alpha\vk{1}\}}}
u^{2/\alpha-\lambda}\int_{x_1}^{x_2}\exp\LT(-\sum_{i=1}^n\frac{\widetilde{f}_i(t)}{\sigma_i^2}\RT)dt
\LT(\prod_{i=1}^n\Psi\LT(\frac{m_{u,i}}{\sigma_{u,i}(0)}\RT)\RT), \ \ u\rw\IF. $$
{\it $\diamond$  \underline{Case $\lambda=2/\alpha$}.} Without loss of generality we assume that $x_1=-\IF$ and $x_2=\IF$.
The  cases  $x_1>-\IF$ and $x_2<\IF$ can be dealt with analogously.
In what follows, we use notation introduced in (\ref{notation1}) and set $\widetilde{I}(u)=I_0(u)\cup I_{-1}(u)$.
Observe that for large $u$
	\BQN
	&&\pi(u)\geq\pk{\exists_{t\in \widetilde{I}(u)} \vk{Z}_{u,-\vn}(t)>\frac{\vk{m}_u}{\vk{\sigma}_u(0)}}\label{proofI0},\\
	&&\pi(u)\leq\pk{ \exists_{t\in  \widetilde{I}(u)} \vk{Z}_{u,+\vn}(t)>\frac{\vk{m}_u}{\vk{\sigma}_u(0)}}+\underset{k\neq -1,0}{\sum_{k=N_1(u)}^{N_2(u)}}\pk{ \exists_{t\in I_{k}(u)} \vk{Z}_{u,+\vn}(t)>\frac{\vk{m}_u}{\vk{\sigma}_u(0)}}.\label{proofI1}
	\EQN
	 Lemma \ref{pan2} yields that
	\BQN\label{proofI2}
	\pk{\exists_{t\in \widetilde{I}(u)} \vk{Z}_{u,\pm\vn}(t)>\frac{\vk{m}_u}{\vk{\sigma}_u(0)}}	\sim\mathcal{P}^{\frac{\widetilde{\vk{f}}}{\vk{\sigma}^2}}
_{\alpha,\frac{\vk{a}}{{\vk{\sigma}}^2}\mathbf{I}_{\{\vk{\alpha}=\alpha\vk{1}\}}}[-S,S]
	\prod_{i=1}^n\Psi\LT(\frac{m_{i,u}}{\sigma_{i,u}(0)}\RT),
	\EQN
	as $\ u\rw\IF, \vn\rw 0$.
	Moreover, in light of \nelem{pan2} and \eqref{FF} we have
	\BQN\label{proofI3}
	&&\underset{k\neq -1,0}{\sum_{k=N_1(u)}^{N_2(u)}}\pk{ \exists_{t\in I_{k}(u)} \vk{Z}_{u,+\vn}(t)>\frac{\vk{m}_u}{\vk{\sigma}_u(0)}}\nonumber\\
	&&\quad\quad \leq\underset{k\neq -1,0}{\sum_{k=N_1(u)}^{N_2(u)}} \pk{ \exists_{t\in I_{0}(u)}\vk{Y}_{+\vn}(t)>\vk{v}_{u,-\vn}(k)}\nonumber\\
	&&\quad\quad\sim\mathcal{H}_{\alpha,(1+\vn)\frac{\vk{a}}{\vk{\sigma}^2}
\mathbf{I}_{\{\vk{\alpha}=\alpha\vk{1}\}}}[0,S]
	\underset{k\neq -1,0}{\sum_{k=N_1(u)}^{N_2(u)}}\prod_{i=1}^n\Psi(v_{u,i,-\vn}(k))\nonumber\\
	&&\quad\quad\sim \mathcal{H}_{\alpha,(1+\vn)\frac{\vk{a}}{\vk{\sigma}^2}\mathbf{I}_{\{\vk{\alpha}=\alpha\vk{1}\}}}[0,S]
	\LT(\prod_{i=1}^n\Psi\LT(\frac{m_{u,i}}{\sigma_{u,i}(0)}\RT)\RT)
	\underset{k\neq -1,0}{\sum_{k=N_1(u)}^{N_2(u)}}\exp\LT(-\sum_{i=1}^n\frac{m_{u,i}^2u^{-2}\inf_{s\in [k,k+1]}\widetilde{f}^\vn_i(S s)}{\sigma^2_{u,i}(0)}\RT)\nonumber\\
	&&\quad\quad \sim\mathcal{H}_{\alpha,(1+\vn)\frac{\vk{a}}{\vk{\sigma}^2}\mathbf{I}_{\{\vk{\alpha}=\alpha\vk{1}\}}}[0,S]
	\LT(\prod_{i=1}^n\Psi\LT(\frac{m_{u,i}}{\sigma_{u,i}(0)}\RT)\RT)
	\underset{k\neq -1,0}{\sum_{k=N_1(u)}^{N_2(u)}}\exp\LT(-\sum_{i=1}^n\frac{\inf_{s\in [k,k+1]}\widetilde{f}^\vn_i(S s)}{\sigma^2_i}\RT)\nonumber\\
	&&\quad\quad \leq \mathbb{C}_4\mathcal{H}_{\alpha,\frac{\vk{a}}{\vk{\sigma}^2}\mathbf{I}_{\{\vk{\alpha}=\alpha\vk{1}\}}}
	\LT(\prod_{i=1}^n\Psi\LT(\frac{m_{u,i}}{\sigma_{u,i}(0)}\RT)\RT)Se^{-\eta\ln S}
	=o\left(\prod_{i=1}^n\Psi\LT(\frac{m_{u,i}}{\sigma_{u,i}(0)}\RT)\right),\ u\rw\IF, \vn\rw 0, S\rw\IF,
	\EQN
where $\eta\in (1,\IF)$ is a constant.
	\COM{Inserting (\ref{proofI2}), (\ref{proofI3}) into (\ref{proofI1}) leads to
		\BQNY
		\lim_{u\rw\IF}\frac{\pi(u)}{\prod_{i=1}^n\Psi\LT(\frac{M_{u,i}}{\sigma_{u,i}(t_u)}\RT)}\leq \mathcal{P}^{\frac{c^2\vk{f}}{\vk{\sigma}^2}}_{\alpha,\vk{a}(\frac{c}{\vk{\sigma}})^2
\mathbf{I}_{\{\vk{\alpha}=\alpha\vk{1}\}}}[-S,S]
		+\mathbb{C}_4\mathcal{H}_{\alpha,\vk{a}(\frac{c}{\vk{\sigma}})^2
\mathbf{I}_{\{\vk{\alpha}=\alpha\vk{1}\}}}Se^{-\mathbb{C}_5 S^{\epsilon}}<\IF,
		\EQNY
		and by (\ref{proofI0}), we have
		\BQNY
		\lim_{u\rw\IF}\frac{\pi(u)}{\prod_{i=1}^n\Psi\LT(\frac{M_{u,i}}{\sigma_{u,i}(t_u)}\RT)}\geq \mathcal{P}^{\frac{c^2\vk{f}}{\vk{\sigma}^2}}_{\alpha,\vk{a}(\frac{c}{\vk{\sigma}})^2
\mathbf{I}_{\{\vk{\alpha}=\alpha\vk{1}\}}}[-S_1,S_1]>0.
		\EQNY}
	Inserting \eqref{proofI2}-\eqref{proofI3} into \eqref{proofI0}-\eqref{proofI1} and letting $S\rw\IF$, we obtain that \COM{$\mathcal{P}^{\frac{c^2\vk{f}}{\vk{\sigma}^2}}_{\alpha,\vk{a}(\frac{c}{\vk{\sigma}})^2
\mathbf{I}_{\{\vk{\alpha}=\alpha\vk{1}\}}}[x_1,x_2]\in(0,\IF)$ and } $$\pi(u)\sim\mathcal{P}^{\frac{\widetilde{\vk{f}}}{\vk{\sigma}^2}}_{\alpha,\frac{\vk{a}}
{\vk{\sigma}^2}\mathbf{I}_{\{\vk{\alpha}=\alpha\vk{1}\}}}
	\prod_{i=1}^n\Psi\LT(\frac{m_{u,i}}{\sigma_{u,i}(0)}\RT),\  u\rw\IF.$$
This establishes the claim.\\
{\it $\diamond$ \underline{Case $\lambda>\frac{2}{\alpha}$}.} Without loss of generality we assume that $x_1=-\IF$ and $x_2=\IF$.
		For any $S>0$, define
		$$J_k(u)=[ku^{-\lambda}S,(k+1)u^{-\lambda}S], k\in\mathbb{Z},\quad \widetilde{J}(u)=J_0(u)\cup J_{-1}(u),$$ $$ K_1(u)=\LT\lfloor \frac{x_1(u)}{S u^{-\lambda}}\RT\rfloor-\mathbb{I}_{\{x_1\leq0\}},\ K_2(u)=\LT\lfloor \frac{x_2(u)}{S u^{-\lambda}}\RT\rfloor+\mathbb{I}_{\{x_2\leq0\}},\quad
		\vk{v}_{u,\pm\vn}(k)=(v_{u,1,+\vn}(k),\ldots,v_{u,n,+\vn}(k)),$$ with
		\BQNY
		v_{u,i,+\vn}(k)=\frac{m_{u,i}}{\sigma_{u,i}(0)}\sup_{s\in J_k(u)}w_{u,i,+\vn}(s),\ \
		v_{u,i,-\vn}(k)=\frac{m_{u,i}}{\sigma_{u,i}(0)}\inf_{s\in J_k(u)}w_{u,i,-\vn}(s).
		\EQNY
		Then for $u$ large enough, we have
		\BQN
		\pi(u)&\geq&\pk{\exists_{t\in \widetilde{J}(u)}\vk{Z}_{u,-\vn}(t)>\frac{\vk{m}_u}{\vk{\sigma}_u(0)}}\label{c31},\\
		\pi(u)&\leq& \pk{ \exists_{t\in \widetilde{J}(u)} \vk{Z}_{u,+\vn}(t)>\frac{\vk{m}_u}{\vk{\sigma}_u(0)}}+\underset{k\neq 0,-1}{\sum_{k=K_1(u)}^{K_2(u)}}\pk{ \exists_{t\in J_{k}(u)} \vk{Z}_{u,+\vn}(t)>\frac{\vk{m}_u}{\vk{\sigma}_u(0)}}.\label{c32}
		\EQN
		It follows from Lemma \ref{pan2} that
		\BQN\label{c33}
		\pk{\exists_{t\in \widetilde{J}(u)} \vk{Z}_{u,\pm\vn}(t)>\frac{\vk{m}_u}{\vk{\sigma}_u(0)}}
		\sim\int_{\R^n}e^{\sum_{i=1}^n w_i}\mathbb{I}_{\LT\{\exists_{t\in[-S,S]}
         -\frac{\widetilde{\vk{f}}(t)}{\vk{\sigma}^2}>\vk{w}\RT\}}d\vk{w}
		\prod_{i=1}^n\Psi\LT(\frac{m_{i,u}}{\sigma_{i,u}(0)}\RT),
		\EQN
		as $u\rw\IF,\ \vn\rw 0$.
		Moreover, similarly to (\ref{proofI3}), we have that
		\BQN\label{c34}
		\underset{k\neq -1,0}{\sum_{k=K_1(u)}^{K_2(u)}}\pk{ \exists_{t\in J_{k}(u)} \vk{Z}_{u,+\vn}(t)>\frac{\vk{m}_u}{\vk{\sigma}_u(0)}}
		&\leq&  \mathbb{C}_6\LT(\prod_{i=1}^n\Psi\LT(\frac{m_{u,i}}{\sigma_{u,i}(0)}\RT)\RT)e^{- \eta \ln S}\nonumber\\
		&=&o\left(\prod_{i=1}^n\Psi\LT(\frac{m_{u,i}}{\sigma_{u,i}(0)}\RT)\right),\ u\rw\IF, S\rw\IF.
		\EQN
	Inserting \eqref{c33}-\eqref{c34} into \eqref{c31}-\eqref{c32}  and letting $S\rw\IF$ and $\epsilon\rw 0$ we derive that
		$$\pi(u)\sim\LT(\int_{\R^n}e^{\sum_{i=1}^n w_i}\mathbb{I}_{\LT\{\exists_{t\in(-\IF,\IF)}
         -\frac{\widetilde{\vk{f}}(t)}{\vk{\sigma}^2}>\vk{w}\RT\}}d\vk{w}\RT)
		\prod_{i=1}^n\Psi\LT(\frac{m_{u,i}}{\sigma_{u,i}(0)}\RT),\  u\rw\IF.$$
This completes the proof. \QED\\

\prooftheo{Thm2} We first focus on the case of  $t_0\in(0,T)$.
Set
$$\E=[-\delta(u), \delta(u)],\quad D(u):=[t_0-\theta,t_0+\theta]\setminus (t_0+\E), $$ where $\theta\in(0,\frac{1}{2})$ is a small constant and
$\delta(u)=\LT(\frac{(\ln u)^q}{u}\RT)^{2/\beta}$   with $q>1$, $\beta=\min_{1\leq i\leq n}\beta^*_i$ and $\beta^*_i=\min\LT(\beta_i, 2\gamma_i\II_{\{c_i\neq0\}}+\IF\II_{\{c_i=0\}}\RT)$.
Then it follows that
\BQNY
\Pi_1(u)\le \pk{\exists_{t\in[0,T]}\LT(\vk{X}(t)+\vk{h}(t)\RT)>u\vk{1}}\le \Pi_1(u)+ \Pi_2(u)+\Pi_3(u),
\EQNY
where
\BQNY
\Pi_1(u)&=&\pk{\exists_{t\in  \E}  \LT(\vk{X}(t_0+t)+\vk{h}(t_0+t)\RT) >u\vk{1}},\quad
\Pi_2(u)=\pk{\exists_{t\in D(u)} \LT(\vk{X}(t)+\vk{h}(t)\RT)>u\vk{1}},\\
\Pi_3(u)&=&\pk{\exists_{t\in [0,T]\setminus [t_0-\theta,t_0+\theta] } \LT(\vk{X}(t)+\vk{h}(t)\RT)>u\vk{1}}.
\EQNY
{\it \underline{Asymptotics of $\Pi_1(u)$}}.
In order to derive  the asymptotics of $\Pi_1(u)$, we check the assumptions in \netheo{PreThm1}.  For this purpose, rewrite
\BQNY
\Pi_1(u)=\pk{\exists_{t\in\E}\vk{X}_u(t)>u\vk{1}}, \quad \text{with} \quad \vk{X}_u(t)=\frac{\vk{X}(t_0+t)}{\vk{1}-{\vk{h}(t_0+t)}/u}.
\EQNY
It follows straightforwardly that $\vk{\sigma}_u(t)=\frac{\vk{\sigma}(t_0+t)}{\vk{1}-\vk{h}(t_0+t)/u}$ satisfies
$\lim_{u\rw\IF}\vk{\sigma}_u(0)=\vk{\sigma}(t_0)>\vk{0}$ implying that {\bf A1} holds.
Next we verify {\bf A2}. Direct calculation shows that
\BQNY
\frac{\sigma_{u,i}(0)}{\sigma_{u,i}(t)}-1=\frac{1}{\sigma_i(t_0+t)} (\sigma_i(t_0)-\sigma_i(t_0+t))+\frac{1}{u-h_i(t_0)}\frac{\sigma_i(t_0)}{\sigma_i(t_0+t)} (h_i(t_0)-h_i(t_0+t)).
\EQNY
 Thus by \eqref{eq:sigt0} and \eqref{eq:gtt0} we have that for all $u$ large
 \BQN \label{eq:Mut0t}
\frac{\sigma_{u,i}(0)}{\sigma_{u,i}(t)}=1+\LT(\frac{b_i}{\sigma_i(t_0)}\abs{t}^{\beta_i} +\frac{c_i}{u-h_i(t_0)}\abs{t}^{\gamma_i}\RT)(1+o(1)),\ t\rw 0.
\EQN
Denote by $\widetilde{f}_i(t)=\frac{{b_i}}{\sigma_i(t_0)}|t|^{\beta_i}\mathbb{I}_{\{\beta_i=\beta^*_i\}}
+c_i|t|^{\gamma_i}\mathbb{I}_{\{\beta^*_i=2\gamma_i\}}
$. Then we have
\BQN\label{corpp1}
\lim_{u\rightarrow\IF}{\sup_{t\in\E}}
\left|\frac{\LT(\frac{\sigma_{u,i}(0)}{\sigma_{u,i}(t)}-1\RT)u^2
-\widetilde{f}_i(u^{2/{\beta^*_i}}t)}
{|\widetilde{f}_i(u^{2/{\beta^*_i}}t)|+1}\right|
=0,
\EQN
which confirms that {\bf A2} is satisfied.
  Apparently, {\bf A3} follows by \eqref{eq:rst}. Thus we conclude that {\bf A1-A3} are satisfied. Also, (\ref{xx}) holds with $x_1=-\IF$ and $x_2=\IF$. Therefore, in light of  \netheo{PreThm1}, we have, as $u\rw\IF$,
\BQN\label{Pi1}
\Pi_1(u)\sim u^{(\frac{2}{\alpha}-\frac{2}{\beta})_{+}}\prod_{i=1}^n\Psi\LT(\frac{u-h_i(t_0)}{\sigma_i(t_0)}\RT)
\LT\{
\begin{array}{ll}
\mathcal{H}_{\alpha,\frac{\vk{a}}{\vk{\sigma}^2(t_0)}\mathbf{I}_{\{\vk{\alpha}
=\alpha\vk{1}\}}}\int_{-\IF}^\IF e^{-\sum_{i=1}^n f_i(x)}dx,&\ \ \text{if}\ \alpha<\beta,\\
\mathcal{P}^{\vk{f}}_{\alpha,\frac{\vk{a}}{\vk{\sigma}^2(t_0)}\mathbf{I}_{\{\vk{\alpha}
=\alpha\vk{1}\}}}(-\IF,\IF), &\ \ \text{if}\ \alpha=\beta,\\
1,&\ \ \text{if}\ \alpha>\beta,
\end{array}
\RT.
\EQN
where
$f_i(t)=\frac{b_i}{\sigma_i^3(t_0)}|t|^{\beta_i} \II_{\{\beta_i=\beta\}}+\frac{c_i}{\sigma_i^2(t_0)}|t|^{\gamma_i} \mathbb{I}_{\{2\gamma_i=\beta\}}, 1\leq i\leq n.$

{\it \underline{Upper bound for $\Pi_2(u)$}}.
%Clearly, we have that $g_u(t):=\sum_{i=1}^n\frac{1}{\sigma^2_{u,i}(t)}$ attains its minimum at the unique point $0$ over $[-t_0,T-t_0]$.
Observe that
\BQN \label{eq:Pitqq}
\Pi_2(u)=\pk{\exists_{t\in D(u)} \LT(\vk{X}(t)+\vk{h}(t)\RT) >u\vk{1}}\leq\pk{\sup_{t\in [-\theta,\theta]\setminus\E} Y_u(t)>u},
\EQN
where
\BQN\label{YY}
Y_u(t)=\sum_{i=1}^nG_{u,i}(t)X_i(t_0+t), \ t\in [-t_0,T-t_0],
\EQN
with
\BQNY
&&G_{u,i}(t):=\LT(\frac{\prod_{j=1,j\neq i}^n\frac{\sigma_j^2(t_0+t)}{(1-h_j(t_0+t)/u)^2}}{A_u(t_0+t)}\RT)\frac{1}{1-h_i(t_0+t)/u},\
\ t\in [-t_0,T-t_0],\\
&&A_u(t)=\sum_{k=1}^n\LT(\prod_{j=1,j\neq k}^n\frac{\sigma_j^2(t)}{(1-h_j(t)/u)^2}\RT),\ \ t\in [0,T].
\EQNY
In order to analyze the variance of $Y_u$, we introduce   $g_u(t)=\sum_{i=1}^n \frac{1}{\sigma_{u,i}^2(t)}$. Using \eqref{eq:Mut0t} we have that
\BQN\label{eq:mm}
g_u(t)-g_u(0)&=&\sum_{i=1}^n\frac{1}{\sigma_{u,i}^2(t)}-\sum_{i=1}^n\frac{1}
{\sigma_{u,i}^2(0)}\nonumber\\
&=&\sum_{i=1}^n\frac{(\sigma_{u,i}(0)-\sigma_{u,i}(t))(\sigma_{u,i}(0)+\sigma_{u,i}(t))}
{\sigma_{u,i}^2(t)\sigma_{u,i}^2(0)}\nonumber\\
&\geq&\mathbb{C}_0\sum_{i=1}^n\frac{1}{\sigma^2(t_0)}\LT(\frac{b_i}{\sigma_i(t_0)}\abs{t}^{\beta_i} +\frac{c_i}{u}\abs{t}^{\gamma_i}\RT)\nonumber\\
&\geq& C\frac{(\ln u)^q}{u^2}
\EQN
holds for all $t\in[-\theta,\theta]\setminus\E$ with a positive constant $C$. Consequently,
\BQNY
\sup_{t\in [-\theta,\theta]\setminus\E}\Var(Y_u(t))=\sup_{t\in [-\theta,\theta]\setminus\E}\LT(\sum_{i=1}^n\frac{(1-h_i(t_0+t)/u)^2}{\sigma_i^2(t_0+t)}\RT)^{-1}
=\sup_{t\in [-\theta,\theta]\setminus\E}\frac{1}{g_u(t)}\leq\frac{1}{g_u(0)+\frac{C(\ln u)^q}{u^2}}.
\EQNY
By (\ref{holderf}) and the fact that in view of (\ref{holderx}),
$$
(\sigma_i(t)-\sigma_i(s))^2\leq \EF{(X_i(t)-X_i(s))^2}\leq \mathbb{C}_1|t-s|^{\mu_1}, \quad s,t\in [0,T],
$$
we have that there exists $\mu_3>0$ such that
\BQNY
\max_{i=1,\dots, n}(G_{u,i}(t)-G_{u,i}(s))^2\leq \mathbb{C}_2|t-s|^{\mu_3},\ s,t\in[0,T],
\EQNY
which together with (\ref{holderx}) implies that
\BQN\label{holdery}
\mathbb{E}\LT(Y_u(t)-Y_u(s)\RT)^2
&=&\mathbb{E}\LT(\sum_{i=1}^nG_{u,i}(t)X_i(t)-\sum_{i=1}^nG_{u,i}(s)X_i(s)\RT)^2\nonumber\\
&=&\sum_{i=1}^n\mathbb{E}\LT(G_{u,i}(t)X_i(t)-G_{u,i}(s)X_i(s)\RT)^2\nonumber\\
&\leq&2\sum_{i=1}^n\sigma_i^2(t)\LT(G_{u,i}(t)-G_{u,i}(s)\RT)^2
+2\sum_{i=1}^n G_{u,i}^2(s)\mathbb{E}\LT(X_i(t)-X_i(s)\RT)^2\nonumber\\
&\leq&\mathbb{C}_3|t-s|^{\mu_4}, \quad s,t\in [0,T]
\EQN
with $\mu_4>0$. Consequently Piterbarg  inequality (Theorem 8.1 in \cite{Pit96}) gives that  %(cf. Theorem 5.3 in \cite{Tabis})
\BQNY
\Pi_2(u)&\leq& \pk{\sup_{t\in [-\theta,\theta]\setminus\E} Y_u(t)>u}\\
&\leq& \mathbb{C}_4 u^{2/\mu_4}\Psi\LT(\sqrt{u^2 g_u(0)+C(\ln u)^q}\RT)\nonumber\\
&=&o\LT(u^{(\frac{2}{\alpha}-\frac{2}{\beta})_{+}}\prod_{i=1}^n
\Psi\LT(\frac{u-h_i(t_0)}{\sigma_{i}(t_0)}\RT)\RT), \quad u\rw\IF.
\EQNY
{\it \underline{Upper bound for $\Pi_3(u)$}}.
Note that there exists $\epsilon\in (0,1)$ such that
\BQNY
\sup_{t\in ([0,T] \setminus  [t_0-\theta,t_0+\theta])} \sigma_i(t)\leq (1-\epsilon)\sigma_i(t_0),\ 1\leq i\leq n.
\EQNY
Thus
\BQNY
\sup_{t\in [0,T]\setminus[-\theta,\theta]}\Var(Y_u(t))
=\left(\inf_{t\in [0,T]\setminus[-\theta,\theta]}g_u(t)\right)^{-1}\leq (1-\epsilon/2)^{-2}\left(\sum_{i=1}^n\frac{1}{\sigma^2_i(t_0)}\right)^{-1},
\EQNY
which together with (\ref{holdery}) and Piterbarg  inequality (Theorem 8.1 in \cite{Pit96}) implies that
\BQNY
 \Pi_3(u)&=& \pk{\exists_{t\in([0,T] \setminus  [t_0-\theta,t_0+\theta])} \LT(\vk{X}(t)+{\vk{h}}(t)\RT) >u\vk{1}}\\
 &\leq&\pk{\sup_{t\in ([0,T] \setminus  [t_0-\theta,t_0+\theta])} Y_u(t)>u}\\
&\leq& \mathbb{C}_5 u^{2/\mu_4}\Psi\left((1-\epsilon/2)\left(\sum_{i=1}^n\frac{1}{\sigma^2_i(t_0)}\right)^{1/2} u\right)\\
 &=&o(\Pi_1(u)),\ \ u\rw\IF.
\EQNY
Therefore, we conclude that
$$\pk{\exists_{t\in[0,T]}\LT(\vk{X}(t)+\vk{h}(t)\RT)>u\vk{1}}\sim  \Pi_1(u), \quad u\rw\IF,$$
which combined with (\ref{Pi1}) establishes the claim.\\
The case of $t_0=0$ ($t_0=T$) can be dealt with using  the same argument as above
with the only difference that one has to substitute $\E$ by  $[0,\delta(u)]$ (or by $[-\delta(u),0]$).

Thus the proof is complete.
\QED

\COM{
%\proofkorr{ex1}
\proofprop{prop} \\
{\it \underline{Ruin probability}}.
We notice that
\begin{align*}
\pk{\exists_{t\in[0,T]}\LT(\vk{B}_{\vk{\alpha}}(t)+\vk{c}t\RT) >u\vk{d}}
=\pk{\exists_{t\in[0,T]}\LT(\frac{1}{\vk{d}}\vk{B}_{\vk{\alpha}}(t)+\frac{\vk{c}t}{\vk{d}}\RT) >u\vk{1}},
\end{align*}
and the variance function $\sigma^2_i(t)$ and correlation function $r_i(s,t)$ of $\frac{B_{\alpha_i}(t)}{d_i}$  satisfy
\begin{align*}
&r_i(s,t)=1-\frac{1}{2T^{\alpha_i}}|t-s|^{\alpha_i}(1+o(1)), s,t\rw T,\\
& \sigma_i(t)=\frac{T^{{\alpha_i}/2}}{d_io}-\frac{\alpha_i}{2d_i}T^{{\alpha_i}/2-1}(T-t)(1+o(1)), t\rw T,
\end{align*}
where $T$ is the unique maximum point of $\sigma_i(t), 1\leq i\leq n$ over $[0,T]$. Moreover,
$$\frac{c_it}{d_i}=\frac{c_iT}{d_i}-\frac{c_i}{d_i}|T-t|, \quad t\rw T.$$
Therefore, in light of \netheo{Thm2}, we have that
\begin{align*}
\pk{\exists_{t\in[0,T]}\LT(\vk{B}_{\vk{\alpha}}(t)+\vk{c}t\RT) >u\vk{d}}\sim u^{(\frac{2}{\alpha}-2)_{+}}\prod_{i=1}^n\Psi\LT(\frac{d_i u-c_iT}{T^{{\alpha_i}/2}}\RT)
\LT\{
\begin{array}{ll}
\mathcal{H}_{\alpha,\vk{e}\mathbf{I}_{\{\vk{\alpha}=\alpha\vk{1}\}}}
\int_{0}^\IF e^{-\sum_{i=1}^n f_i(t)}dt,&\ \ \text{if}\ \alpha<1,\\
\mathcal{P}^{\vk{f}}_{\alpha,\vk{e}
\mathbf{I}_{\{\vk{\alpha}=\alpha\vk{1}\}}}[0,\IF), &\ \ \text{if}\ \alpha=1,\\
1,&\ \ \text{if}\ \alpha>1,
\end{array}
\RT.
\end{align*}
where $\alpha=\min_{1\leq i\leq n}\alpha_i$, \bl{$\vk{e}=(e_1,\dots, e_n)$ with $e_i=\frac{d_i^2}{2T^{2\alpha_i}}$ }and $f_i(t)=\frac{\alpha_id_i^2}{2T^{\alpha_i+1}}|t|$.\\
{\it \underline{Ruin time}}.
By definition,
\BQN\label{ratio}
\pk{(T-\tau_u)u^2\leq x\Big| \tau_u\leq T}=\frac{\pk{\exists_{t\in[T-xu^{-2},T]} \LT(\vk{B}_{\vk{\alpha}}(t)+\vk{c}t\RT) >u\vk{d}}}{\pk{\exists_{t\in[0,T]} \LT(\vk{B}_{\vk{\alpha}}(t)+\vk{c}t\RT) >u\vk{d}}}.
\EQN
For all $u>0$ large,
\BQNY
\pk{\exists_{t\in[T-xu^{-2},T]} \LT(\vk{B}_{\vk{\alpha}}(t)+\vk{c}t\RT) >u\vk{d}}
&=&\pk{\exists_{t\in[-xu^{-2},0]} \LT(\vk{B}_{\vk{\alpha}}(T+t)+\vk{c}(T+t)\RT) >u\vk{d}}\\
&=&\pk{\exists_{t\in[-xu^{-2},0]}\frac{ \vk{B}_{\vk{\alpha}}(T+t)}{\vk{d}-\vk{c}(T+t)/(u\vk{1})} >u\vk{1}}
\EQNY
Denote by $\vk{X}_u(t)=\frac{ \vk{B}_{\vk{\alpha}}(T+t)}{\vk{d}-\vk{c}(T+t)/(u\vk{1})}$.
Let $\sigma_{u,i}(t)$ and $r_{u,i}(s,t)$ denote the variance and correlation functions of $X_{u,i}(t)$, respectively. Next we check the conditions of \netheo{PreThm1}. {\bf A1} holds straightforwardly. Direct calculation gives that
 $$
 \lim_{u\rightarrow\IF}{\sup_{t\in [-xu^{-2},0]}}\LT|\frac{\LT(\frac{\sigma_{u,i}(0)}{\sigma_{u,i}(t)}-1\RT)u^2-f
 (u^2t)}
 {\LT|f(u^2t)\RT|+1}\RT|=0,
$$
 where $f(t)=\frac{\alpha}{2T}|t|$, and
$$
\lim_{u\rightarrow\IF}\sup_{s,t\in [-xu^{-2},0]}\left|\frac{1-r_{u,i}(t,s)}{2^{-1}T^{-\alpha_i}|t-s|^{\alpha_i}}-1\right|=0.
$$
These confirm {\bf A2-A3} hold.
Moreover, it easy to check that (\ref{xx}) is satisfied with $x_1=-x$ and $x_2=0$. Consequently, in light of  \netheo{PreThm1}, we have that, as $u\rw\IF$,
\BQNY
\pk{\exists_{t\in[T-xu^{-2},T]} \LT(\vk{B}_{\vk{\alpha}}(t)+\vk{c}t\RT) >u\vk{d}}\sim u^{(\frac{2}{\alpha}-2)_{+}}\prod_{i=1}^n\Psi\LT(\frac{d_i u-c_iT}{T^{{\alpha_i}/2}}\RT)
\LT\{
\begin{array}{ll}
\mathcal{H}_{\alpha, \vk{e}\mathbf{I}_{\{\vk{\alpha}=\alpha\vk{1}\}}}
\int_{-x}^0 e^{-\sum_{i=1}^n f_i(t)}dt,&\ \ \text{if}\ \alpha<1,\\
\mathcal{P}^{\vk{f}}_{\alpha, \vk{e}\mathbf{I}_{\{\vk{\alpha}
=\alpha\vk{1}\}}}[-x,0], &\ \ \text{if}\ \alpha=1,\\
1,&\ \ \text{if}\ \alpha>1,
\end{array}
\RT.
\EQNY
with the same $\alpha, \vk{e}$ and $\vk{f}$ as in the proof of (\ref{e.1}).
Inserting the above asymptotics into (\ref{ratio}), we establish the claim.
\QED
\COM{
\proofkorr{ex2}
Set $\vk{X}(t)=(X_1(t),\cdots,X_n(t))$ with $X_i(t)=\frac{B^i_\alpha(t)}{d-c t}, \ 1\leq i\leq n$ and $\mu=u^{1-\alpha/2}$. Then we have for $u>0$
\begin{align*}
\pk{\exists_{t\in[0,\IF)}\LT(\vk{B}_\alpha(t)+\vk{c}t\RT)>u\vk{d}}
=\pk{\exists_{t\in[0,\IF)}\vk{X}(t)>\mu\vk{1}}.
\end{align*}
$X_i(t)$'s variance function  $\sigma^2(t)=\frac{t^{\alpha}}{(d-c t)^2}$ attain its maximum over $[0,\IF)$ at point $t_0=\frac{d\alpha}{c(\alpha-2)}.$ Denoting  $\sigma_0=\sigma(t_0)=\frac{t_0^{\alpha/2}}{d-ct_0}$, we have
\BQN
&&\sigma(t)=\sigma_0-\frac{\alpha(2-\alpha)t_0^{\alpha/2-2}}{8(d-ct_0 )}(t-t_0)^2+o((t-t_0)^2), \ t \rw t_0,\label{EX2v}\\
%&\frac{\sigma_i(t_0)}{\sigma_i(t)}-1=\frac{2\alpha d}{8(ct_0^3+dt_0^2)}(t-t_0)^2+o((t-t_0)^2), \ 1\leq i\leq n.\\
&&r(s,t)=1-\frac{1}{2t_0^{\alpha}}|t-s|^{\alpha}(1+o(1)),\  s,t\rw t_0.\label{EX2r}.
\EQN

We have for $K>t_0$
\BQNY
&&\pk{\exists_{t\in[0,\IF)}\vk{X}(t)>\mu\vk{1}}
\geq\pk{\exists_{t\in[0,K]}\vk{X}(t)>\mu\vk{1}},\\
&&\pk{\exists_{t\in[0,\IF)}\vk{X}(t)>\mu\vk{1}}
\leq\pk{\exists_{t\in[0,K]}\vk{X}(t)>\mu\vk{1}}
+\pk{\exists_{t\in[K,\IF)}\vk{X}(t)>\mu\vk{1}}.
\EQNY
By \eqref{EX2v}, \eqref{EX2r} and \netheo{Thm2} we have as $u\rw\IF$
\BQNY
\pk{\exists_{t\in[0,K]}\vk{X}(t)>\mu\vk{1}}
&\sim& \mu^{\frac{2}{\alpha}-2}\LT(\Psi\LT(\frac{\mu}{\sigma_0}\RT)\RT)^n
\mathcal{H}_{\alpha,a\vk{1}}\int_{-\IF}^\IF e^{-n f(t)}dt
=\frac{\sqrt{2\pi}}{q}\mathcal{H}_{\alpha,a\vk{1}}
\mu^{\frac{2}{\alpha}-2}\LT(\Psi\LT(\frac{\mu}{\sigma_0}\RT)\RT)^n,
\EQNY
where $a=\frac{d^2(d-ct_0)^2}{2t_0^{2\alpha}}$, and $f(t)=\frac{q^2}{2n}
t^2$.
By Borell-TIS inequality and
\BQNY
V:=\sup_{t\in[K,\IF)}\frac{1}{n}\sigma^2(t)<\frac{\sigma^2_0 }{n},
\EQNY
we know
\BQNY
\pk{\exists_{t\in[K,\IF)}\vk{X}(t)>\mu\vk{1}}&\leq&
 \pk{\exists_{t\in[K,\IF)}\frac{\sum_{i=1}^nX_i(t)}{n}>\mu}\\
 &\leq& \exp\LT(-\frac{(\mu-\mathbb{Q})^2}{2V}\RT)\\
 &=&o\LT(\Psi^n\LT(\frac{\mu}{\sigma_0}\RT)\RT),
\EQNY
where
\BQNY
\mathbb{Q}=\mathbb{E}\LT\{\sup_{t\in[0,\IF)}\frac{\sum_{i=1}^nX_i(t)}{n}\RT\}<\IF.
\EQNY
Thus the claim follows.
\QED

\proofkorr{ex2ruin} Let $\vk{X}(t)$ and $\mu$ be the same as in the proof of \nekorr{ex2}.
By definition, for $x\in \R$
\BQNY
\pk{q(\tau_u-t_0 u)\leq x\Big| \tau_u< \IF}=\frac{\pk{\exists_{t\in[0,x/q+t_0u]}\LT(\vk{B}_\alpha(t)+\vk{c}t\RT)>u\vk{d}}}
{\pk{\exists_{t\in[0,\IF)}\LT(\vk{B}_\alpha(t)+\vk{c}t\RT)>u\vk{d}}}
=\frac{\pk{\exists_{t\in[0,\frac{x}{q u}+t_0]}\vk{X}(t)>\mu\vk{1}}}
{\pk{\exists_{t\in[0,\IF)} \vk{X}(t)>\mu\vk{1}}}.
\EQNY
 As in the proof of \netheo{Thm2}, by \netheo{PreThm1} we obtain as $u\rw\IF$
\BQNY
\pk{\exists_{t\in[0,\frac{x}{q u}+t_0]}\vk{X}(t)>\mu\vk{1}}\sim \mu^{\frac{2}{\alpha}-2}\LT(\Psi\LT(\frac{\mu}{\sigma_0}\RT)\RT)^n
\mathcal{H}_{\alpha,a\vk{1}}\int_{-\IF}^{x/q} e^{-n f(t)}dt,
\EQNY
where $a=\frac{d^2(d-ct_0)^2}{2t_0^{2\alpha}}$ and $f(t)=\frac{q^2}{2n}
t^2$.
Thus by \nekorr{ex2}, we have as $u\rw\IF$
\BQNY
\pk{q(\tau_u-t_0 u)\leq x\Big| \tau_u< \IF}\sim\frac{\int_{-\IF}^{x/q} e^{-n f(t)}dt}{\int_{-\IF}^{\IF} e^{-n f(t)}dt}.
\EQNY
\QED
}

\COM{		
{\bf  Proof of Example \ref{ex1}:}
we have
\begin{align*}
\pk{\exists_{t\in[0,T]}\LT(\vk{B}_{\vk{\alpha}}(t)+\vk{c}t\RT) >u\vk{d}}
=\pk{\exists_{t\in[0,T]}\LT(\frac{1}{\vk{d}}\vk{B}_{\vk{\alpha}}(t)+\frac{\vk{c}t}{\vk{d}}\RT) >u\vk{1}}.
\end{align*}
Let $\sigma^2_i(t)$ and $r_i(t)$ be the variance function and correlation function of $\frac{1}{d_i}B_{\alpha_i}(T)$ respectively. Then $\sigma_i(t), 1\leq i\leq n$ attain its unique maximizer at point $t=T$ over $[0,T]$.
Further,
\begin{align*}
&r_i(s,t)=1-\frac{1}{2T^{\alpha_i}}|t-s|^{\alpha_i}(1+o(1)), s,t\rw T,\\
& \sigma_i(t)=\frac{T^{{\alpha_i}/2}}{d_i}-\frac{\alpha_i}{2d_i}T^{{\alpha_i}/2-1}(T-t)(1+o(1)), t\rw T.
\end{align*}
By \netheo{Thmm2}, the results follow.
\QED}

	%%%%%%%%%%%%%%%%%%%%%%%%%%%%%%%%%%%%%%%%%%%%%%%5
	\COM{
		\begin{remark}1
			In \netheo{PreThm1}, we are not able to derive the exact asymptotics when $\alpha>\beta^*$, by using the classical Pickands-Piterbarg technique (see, e.g., \cite{Pit96}). Below we derive some bounds for this case.
			\BQN
			\pk{\sup_{t\in\E} X_g(t_u+t)>u}\geq\pk{  X_g(t_u) >u}= \Psi(M_u^*).
			\EQN
			For any $\vn >0$ we have  $\delta_u\le \vn  u^{-\frac{2}{\alpha}}$ for all large $u$.  Thus, using similar arguments as in the proof of \netheo{PreThm1} we have
			\BQNY
			\pk{\sup_{t\in\E} X_g(t_u+t)>u} \leq \pk{\sup_{t\in [-d_2u^{-\alpha/2}, \vn u^{-\alpha/2}]} X_g(t_u+t)>u}\le \mathcal{H}_\alpha[0,a^{\frac{1}{\alpha}}(d_2+\vn)] \Psi(M_u^*)\oo.
			\EQNY
			Consequently, letting $\vn \to 0$ we obtain
			$$
			1\le \lim_{u\to\IF}\frac{\pk{\sup_{t\in\E}X_g(t_u+t)>u}}{\Psi(M_u^*)}\le \mathcal{H}_\alpha[0,a^{\frac{1}{\alpha}} d_2 ].
			$$
			
		\end{remark}
	}
	%%%%%%%%%%%%%%%%%%%%%%%%%%%%%%%%%%%%%%%%%%%%%%%%%%%%%%%%%%%%%%%%%%%%%%55
\COM{
\BEL\label{lem1}
Let $f\in C_0^*(E)$ and $E$ is a closed interval in $\R$, then $\int_{E} e^{-cf(t)}dt<\IF$ holds for any constant $c>0$.
\EEL
\prooflem{lem1} Since by $f\in C_0^*(E)$ the function $f$ is locally bounded and integrable, it suffice to show that
$\int_b^\IF e^{-cf(t)}dt<\IF$ and $\int_{-\IF}^{-b} e^{-cf(t)}dt<\IF$  for some $b$ large enough. By the assumption, we can find $a\in(0,\IF)$, such that
$cf(t)>t^{\epsilon_1}$ holds for any $t>a$. Hence
\BQNY
\int_a^\IF e^{-c_1f(t)}dt\leq
\int_{0}^\IF e^{-t^{\epsilon_1}}dt=\frac{1}{\epsilon_1}\Gamma(\frac{1}{\epsilon_1})
<\IF.
\EQNY
Similarly, we have $\int_{-\IF}^{-b} e^{-cf(t)}dt<\IF$.
Thus we established the claim. \QED}
\COM{Next we give the tail probability of vector-valued Gaussian processes with trend. Let ${\bf X}(t),t\in[0,T]$ be a centered vector-valued Gaussian processes with independent coordinates $X_{i}$'s which have continuous trajectories,  zeros means, variance functions $\sigma^2_{i}(\cdot), 1\leq i\leq n$ and correlation functions $r_{i}(\cdot,\cdot), 1\leq i\leq n$.\\
Consider the continuous trend function $\vk{h}(t)$ where $h_i(t), 1\leq i\leq n$ bounded measurable functions. In this section, we present the main results concerning the asymptotics of
\BQN\label{eq:Xgt1}
\pk{\exists_{t\in[0,T]}\LT(\vk{X}(t)+\vk{h}(t)\RT)>u\vk{1}}, \ \ \ u\to\IF.
\EQN
Further, we suppose:\\
\textbf{A1'}:The following generalized variance function
\BQNY
\widetilde{g}_u(t):=\sum_{i=1}^{n}\LT(\frac{1-{h_i(t)}/{u}}{\sigma_i(t)}\RT)^2
\EQNY
attains its minimum over $[0,T]$  at the unique point $t_u$ and
$$
\lim_{u\to\IF} t_u=t_0\in[0,T].
$$
\textbf{A4}: For all u large enough,
$$
\inf_{t\in[0,T]\setminus(t_u+\E)} (\widetilde{g}_u(t)-\widetilde{g}_u(t_u))>\frac{p(\ln u)^q}{u^2}
$$
holds for some constants $q>1, p>0$.

\textbf{A5}: For some positive constants $G, \varsigma$,
\BQNY
\sum_{i=1}^n\E{(X_i(t)-X_i(s))^2}\leq G|t-s|^\varsigma,\ \ \sum_{i=1}^n(h_i(t)-h_i(s))^2\leq G|t-s|^\varsigma,
\EQNY
holds for all $s,t\in[0,T]$.\\
In the following theorem, denote $\delta_u=\LT(\frac{(\ln u)^q}{u}\RT)^\lambda$ with $q>1$ in {\bf A4} and
\BQN\label{Du}
\E=
\LT\{
\begin{array}{ll}
{[0,\delta_u]} & \  \text{if} \  \ t_u\equiv 0,\\
{[-t_u,\delta_u]}, & \  \text{if} \  \ t_u\sim du^{-\nu}\ \text{and}\ \nu\geq\lambda,\\
{[-\delta_u,\delta_u]}, & \  \text{if}\ \ t_u\sim du^{-\nu}\ \text{or}\ T-t_u\sim du^{-\nu},\ \nu\leq\lambda\ \text{or}\ t_0\in(0,T),\\
{[-\delta_u,T-t_u]}, & \  \text{if} \  \ T-t_u\sim du^{-\nu}\ \text{and}\ \nu\geq\lambda,\\
{[-\delta_u,0]} & \  \text{if} \  \ t_u\equiv T,
\end{array}
\RT.
\EQN
for $\nu\in(0,\IF)$ and some constant $d>0$.

\BT\label{MainThm1}
For ${\bf X}(t),0\leq t\leq T$ and $\vk{h}(t)$ above, suppose that \textbf{A1'},\textbf{A2}-\textbf{A5} are satisfied with $\sigma_{u,i}(\cdot)=\frac{\sigma_i(t+t_u)}{1-{h_i(t+t_u)}/{u}}$, $r_{u,i}(s,t)=r_i(t_u+s,t_u+t)$, , $\E$ in (\ref{Du}) and $[x_1,x_2]=\lim_{u\rw\IF} u^\lambda\E$ is well defined. Set $\lambda=\max_{1\leq i\leq n}\lambda_{i}$, $\alpha=\min_{1\leq i\leq n}\alpha_i$,  and $\widehat{\vk{f}}(t)=(\widehat{f}_1(t),\cdots,\widehat{f}_n(t))$ with $\widehat{f}_i(t)= f_{i}\LT(\II_{\{\lambda_{i,1}=\lambda\}}t\RT)$, $\widehat{f}_i(t)\in C_0^*([x_1,x_2])$. Further if
$$\lim_{|t|\rw\IF,t\in [x_1,x_2]}\LT(\sum_{i=1}^n\frac{f_i(t)}{\sigma_i^2(t_0)}\RT)\big/{\ln |t|}={\epsilon}\in(1,\IF),$$ then we have
\BQNY\label{ei2}
\pk{\exists_{t\in[0,T]}\LT(\vk{X}(t)+\vk{h}(t)\RT)>u\vk{1}}\sim u^{(\frac{2}{\alpha}-\lambda)_{+}}\prod_{i=1}^n\Psi\LT(\frac{u-h_i(t_u)}{\sigma_{i}(t_u)}\RT)
\left\{
\begin{array}{ll}
	\mathcal{H}_{\alpha,\frac{\vk{a}}{\vk{\sigma}^2(t_0)}\mathbf{I}_{\{\vk{\alpha}=\alpha\vk{1}\}}}
	\int_{x_1}^{x_2}e^{-\sum_{i=1}^n\frac{\widehat{f}_i(t)}{\sigma_i^2(t_0)}}dt, &\hbox{if} \ \ \lambda<2/\alpha,\\
	\mathcal{P}^{\frac{\widehat{\vk{f}}}{\vk{\sigma}^2(t_0)}}_{\alpha,\frac{\vk{a}}
{\vk{\sigma}^2(t_0)}\mathbf{I}_{\{\vk{\alpha}=\alpha\vk{1}\}}}[x_1,x_2],& \hbox{if} \ \ \lambda=2/\alpha,\\
	1,&  \hbox{if} \ \ \lambda>2/\alpha.
\end{array}
\right.
\EQNY
\ET
\COM{
	\BQNY
	d_1^*:=
	\left\{
	\begin{array}{ll}
		0, & \text{if  \textbf{C1} holds},\\
		d, & \text{if \textbf{C2} holds},\\
		\IF, & \text{if \textbf{C3} holds},
	\end{array}
	\right.
	\EQNY
	where\\
	\textbf{C1}: $t_u\sim du^{-\nu}$, or $T-t_u\sim du^{-\nu}$, $0<\lambda<\nu\leq\IF$;\\
	\textbf{C2}: $t_u\sim du^{-\nu}$, or $T-t_u\sim du^{-\nu}$, $0<\nu=\lambda<\IF$;\\
	\textbf{C3}: $t_u\sim du^{-\nu}$, or $T-t_u\sim du^{-\nu}$, $0\leq\nu<\lambda<\IF$ or $t_0\in(0,T)$,\\
	with some positive constant $d$ and $\nu\in(0,\IF]$, $\delta(u)=du^{-\nu}$ if \textbf{C1} and \textbf{C2} hold, $\delta(u)=\delta_u$ if \textbf{C3} holds.}
\COM{\begin{remark}
	i) In \netheo{MainThm1}, $\nu=\IF$ means $t_u\equiv t_0\in\{0,T\}$.
When $t_u\equiv t_0\in[0,T]$, the assumptions \textbf{A4} and \textbf{A5} can be weakened significantly as: \\
		\textbf{A4'}: For any $t\in \Delta_\theta\setminus (t_0+\E)$
		$$
		M_u(t)\geq  M_u(s)
		$$
		holds for any $s\in\E$ and $\theta>0$ when $u$ is large enough.\\
		\textbf{A5'}: For some positive $G, \varsigma$ ($\overline{X}(t)=\frac{X(t)}{\sigma(t)}$),
		$$\E{( \overline{X}(t)-\overline{X}(s))^2}\leq G|t-s|^\varsigma$$
		holds for all $s,t\in\Delta_\theta$ and the same $\theta$ in \textbf{A4'}.\\
		Where $\Delta_\theta=[t_0-\theta,t_0+\theta]\cap[0,T]$ for some $\theta>0$.
\end{remark}}

\prooftheo{MainThm1}
Clearly
\BQNY
\Pi_1(u)\le \pk{\exists_{t\in[0,T]}\LT(\vk{X}(t)+\vk{h}(t)\RT)>u\vk{1}}\le \Pi_1(u)+ \Pi_2(u),
\EQNY
where
\BQNY
\Pi_1(u):=\pk{\exists_{t\in  \E}  \LT(\vk{X}(t_u+t)+\vk{h}(t_u+t)\RT) >u\vk{1}},\
\Pi_2(u):=\pk{\exists_{t\in D(u)} \LT(\vk{X}(t)+\vk{h}(t)\RT)>u\vk{1}}, \ D(u):=[0,T]\setminus (t_u+\E).
\EQNY
Next, we give a sharp upper bound  for $\Pi_2(u)$ which will finally imply that
\BQN\label{eq:Pi2.1}
\Pi_2(u)=o(\Pi_1(u)),\ \ \ u\to\IF.
\EQN
By {\bf A4}, we have that when $u$ is large enough
\BQNY
\inf_{t\in D(u)}\widetilde{g}_u(t)\geq \widetilde{g}_u(t_u)+\frac{p(\ln u)^q}{u^2}
\EQNY
holds for all $t\in D(u)$.
\BQN \label{eq:Pitqq}
\Pi_2(u)=\pk{\exists_{t\in D(u)} \LT(\vk{X}(t)+\vk{h}(t)\RT) >u\vk{1}}\leq\pk{\sup_{t\in D(u)} Y_u(t)>u},
\EQN
where
\BQN\label{YY}
Y_u(t)=\sum_{i=1}^nG_{u,i}(t)X_i(t), \ t\in D(u),
\EQN
with
\BQNY
G_{u,i}(t):=\LT(\frac{\prod_{j=1,j\neq i}^n\frac{\sigma_j^2(t)}{(1-h_j(t)/u)^2}}{A_u(t)}\RT)\frac{1}{1-h_i(t)/u},\
A_u(t)=\sum_{k=1}^n\LT(\prod_{j=1,j\neq k}^n\frac{\sigma_j^2(t)}{(1-h_j(t)/u)^2}\RT),\ t\in D(u).
\EQNY
Since
\BQNY
\Var(Y_u(t))=\LT(\sum_{i=1}^n\frac{(1-h_i(t)/u)^2}{\sigma_i^2(t)}\RT)^{-1}
=\frac{1}{\widetilde{g}_u(t)}\leq\frac{1}{\widetilde{g}_u(t_u)+\frac{p(\ln u)^q}{u^2}},
\EQNY
for all $t\in D(u)$ and by {\bf A5} and H\"older inequality, we know for $u$ large enough
\BQNY
(\sigma_i(t)-\sigma_i(s))^2&=&\sigma_i^2(t)+\sigma_i^2(s)-2\sigma_i(t)\sigma_i(s)\\
&\leq&\sigma_i^2(t)+\sigma_i^2(s)-2\E{X_i(s)X_i(t)} \\
&=&\E{(X_i(t)-X_i(s))^2}\\
&\leq& G|t-s|^\varsigma,\ \  s,t\in[0,T],
\EQNY
and
\BQNY
(G_{u,i}(t)-G_{u,i}(s))^2\leq \mathbb{C}_1|t-s|^\varsigma,\ s,t\in[0,T].
\EQNY
Further,
\BQNY
\mathbb{E}\LT(Y_u(t)-Y_u(s)\RT)^2
&=&\mathbb{E}\LT(\sum_{i=1}^nG_{u,i}(t)X_i(t)-\sum_{i=1}^nG_{u,i}(s)X_i(s)\RT)^2\\
&=&\sum_{i=1}^n\mathbb{E}\LT(G_{u,i}(t)X_i(t)-G_{u,i}(s)X_i(s)\RT)^2\\
&\leq&2\sum_{i=1}^n\mathbb{E}\LT(G_{u,i}(t)X_i(t)-G_{u,i}(s)X_i(t)\RT)^2
+2\sum_{i=1}^n\mathbb{E}\LT(G_{u,i}(s)X_i(t)-G_{u,i}(s)X_i(s)\RT)^2\\
&=&2\mathbb{C}_2\sum_{i=1}^n\LT(G_{u,i}(t)-G_{u,i}(s)\RT)^2
+2\sum_{i=1}^n G_{u,i}^2(s)\mathbb{E}\LT(X_i(t)-X_i(s)\RT)^2\\
&\leq&\mathbb{C}_3|t-s|^{\varsigma},
\EQNY
for all $s,t\in[0,T]$, then by Piterbarg's inequality %(cf. Theorem 5.3 in \cite{Tabis})
\BQN
\pk{\sup_{t\in D(u)} Y_u(t)>u}
&\leq& \mathbb{C}_4T u^{2/\varsigma}\Psi\LT(\sqrt{u^2 \widetilde{g}_u(t_u)+p(\ln u)^q}\RT)\nonumber\\
&=&o\LT(\prod_{i=1}^n\Psi\LT(\frac{u-h_i(t_u)}{\sigma_{i}(t_u)}\RT)\RT),
\EQN
as $u\rw\IF$, where we use the fact that $u^2 \widetilde{g}_u(t_u)\sim \sum_{i=1}^n\frac{u^2}{\sigma_i(t_0)}$, $u\rw\IF$.\\
Next we present the asymptotics of $\Pi_1(u)$. Since
\BQNY
\Pi_1(u)=\pk{\exists_{t\in\E}\frac{\vk{X}(t_u+t)}{\vk{1}-\frac{\vk{h}(t_u+t)}{u\vk{1}}}>u\vk{1}},
\EQNY
then if we set $\vk{X}_u(t)=\frac{\vk{X}(t_u+t)}{\vk{1}-\frac{\vk{h}(t_u+t)}{u\vk{1}}}$, $\vk{\sigma}_u(t)=\frac{\vk{\sigma}(t_u+t)}{\vk{1}-\vk{h}(t_u+t)/(u\vk{1})}$ is the standard deviation with
$\lim_{u\rw\IF}\vk{\sigma}_u(0)=\vk{\sigma}(t_0)$ which satisfy assumption {\bf A2}, and {\bf A1'} equals to {\bf A1} with $g_u(t)=\widetilde{g}_u(t_u+t)$ and $E_u=[-t_u,T-t_u]$.\\
Then by \netheo{PreThm1}, we get
\BQNY
\Pi_1(u)\sim\prod_{i=1}^n\Psi\LT(\frac{u-h_i(t_u)}{\sigma_{i}(t_u)}\RT)\left\{
\begin{array}{ll}
	\mathcal{H}_{\alpha,\frac{\vk{a}}{\vk{\sigma}^2(t_0)}\mathbf{I}_{\{\vk{\alpha}=\alpha\vk{1}\}}}
	u^{\frac{2}{\alpha}-\lambda}\int_{x_1}^{x_2}e^{-\sum_{i=1}^n\frac{\widehat{f}_i(t)}
{\sigma_i^2(t_0)}}dt, &\hbox{if} \ \ \lambda<2/\alpha,\\
	\mathcal{P}^{\frac{\widehat{\vk{f}}}{\vk{\sigma}^2(t_0)}}
_{\alpha,\frac{\vk{a}}{\vk{\sigma}^2(t_0)}\mathbf{I}_{\{\vk{\alpha}=\alpha\vk{1}\}}}[x_1,x_2],& \hbox{if} \ \ \lambda=2/\alpha,\\
	1,&  \hbox{if} \ \ \lambda>2/\alpha,
\end{array}
\right.
\EQNY
as $u\rw\IF$, which further  implies that \eqref{eq:Pi2.1} holds.\\
Thus the proof is complete.
\QED}
}

\prooftheo{Thm3} i) We provide the proof only for case $t_0\in(0,T)$, since cases $t_0=0$ and $t_0=T$ can be established analogously.
Let $\E=[-\delta(u),\delta(u)]$, where $\delta(u)=\LT(\frac{(\ln u)^q }{u}\RT)^{1/\gamma}$ with $q>1$. It follows that
\BQNY
\Pi(u)\le \pk{\exists_{t\in[0,T]}  \LT(\vk{X}(t)+{\vk{h}}(t)\RT)  >u\vk{1}}\le \Pi(u)+ \Pi_1(u),
\EQNY
where
\BQNY
 \Pi(u)=\pk{\exists_{t\in \E} \LT(\vk{X}(t_0+t)+{\vk{h}}(t_0+t)\RT)  >u\vk{1}}, \quad \Pi_1(u)=\pk{\exists_{t\in [0,T] \setminus  (t_0+\E) }   \LT(\vk{X}(t)+{\vk{h}}(t)\RT) >u\vk{1}}.
\EQNY
In order to derive the asymptotics of $\Pi(u)$ we apply \netheo{PreThm1} by checking conditions {\bf A1-A3}.
Set $ \sigma_{u,i}(t)=\frac{1}{1-h_i(t_0+t)/u}$ and then $\lim_{u\rw\IF}\sigma_{u,i}(0)=1$, which indicates that {\bf A1} holds.  By the fact that
\BQNY
\frac{\sigma_{u,i}(0)}{\sigma_{u,i}(t)}-1=\frac{ h_i(t_0)-h_i(t_0+t)}{u-h_i(t_0)},
\EQNY
and \eqref{eq:gtt0}, we have
\BQNY\label{corpp2}
%\lim_{u\rightarrow\IF}\underset{t\not= 0}{\sup_{t\in\E}}\left|\frac{\frac{\sigma_{u,i}(0)}{\sigma_{u,i}(t)}-1}{c_iu^{-1}|t|^{\gamma_i}}-1\right|
\lim_{u\rightarrow\IF}\underset{t\not= 0}{\sup_{t\in\E}}\left|\frac{\LT(\frac{\sigma_{u,i}(0)}{\sigma_{u,i}(t)}-1\RT)u^{2}
-c_i|u^{\frac{1}{\gamma_i}}t|^{\gamma_i}}{c_i|u^{\frac{1}{\gamma_i}}t|^{\gamma_i}+1}\right|
=0.
\EQNY
This confirms that {\bf A2} is satisfied.
Moreover,  \eqref{stationaryR0} implies that
\BQNY
\lim_{u\rightarrow\IF}\underset{t\not=s}{\sup_{  t\in\E,  s\in\E}}\left|\frac{1-r_i(t_0+t, t_0+s)}{ a_i(t_0)|t-s|^{\alpha_i}}-1\right|= 0,
\EQNY
which means that {\bf A3} holds.  Also, we have that (\ref{xx}) holds with $x_1=-\IF$ and $x_2=\IF$.
Therefore, by \netheo{PreThm1}
\BQNY
\Pi(u)\sim u^{(\frac{2}{\alpha}-\frac{1}{\gamma})_{+}}\prod_{i=1}^n\Psi\LT(u-h_{m,i}\RT)
\LT\{
\begin{array}{ll}
\mathcal{H}_{\alpha,\vk{a}_0\mathbf{I}_{\{\vk{\alpha}=\alpha\vk{1}\}}}\int_{-\IF}^\IF e^{-\sum_{i=1}^n f_i(x)}dx,&\ \ \text{if}\ \alpha<2\gamma,\\
\mathcal{P}^{\vk{f}}_{\alpha,\vk{a}_0\mathbf{I}_{\{\vk{\alpha}=\alpha\vk{1}\}}}(-\IF,\IF), &\ \ \text{if}\ \alpha=2\gamma,\\
1,&\ \ \text{if}\ \alpha>2\gamma,
\end{array}
\RT.
\EQNY
where  $\gamma=\min_{1\leq i\leq n}( \gamma_i\II_{\{c_i\neq0\}}+\IF\II_{\{c_i=0\}})$,
$f_i(t)=c_i|t|^\gamma \mathbb{I}_{\{\gamma_i=\gamma\}}, 1\leq i\leq n$.
Next we show that
$
\Pi_1(u)=o(\Pi(u)), u\rw\IF.
$
Observe that
\BQNY
\Pi_1(u)=\pk{\exists_{t\in [0,T] \setminus  (t_0+\E)} \LT(\vk{X}(t)+\vk{h}(t)\RT) >u\vk{1}}\leq\pk{\sup_{t\in [0,T] \setminus  (t_0+\E)} Y_u(t)>u},
\EQNY
where
\BQN\label{Yu}
Y_u(t)=\sum_{i=1}^nG_{u,i}(t)X_i(t_0+t), \ t\in [-t_0,T-t_0],
\EQN
with
\BQNY
&&G_{u,i}(t):=\LT(\frac{\prod_{j=1,j\neq i}^n\frac{1}{(1-h_j(t_0+t)/u)^2}}{A_u(t_0+t)}\RT)\frac{1}{1-h_i(t_0+t)/u},\
\ t\in [-t_0,T-t_0],\\
&&A_u(t)=\sum_{k=1}^n\LT(\prod_{j=1,j\neq k}^n\frac{1}{(1-h_j(t)/u)^2}\RT),\ \ t\in [0,T].
\EQNY
Let
\BQN\label{gt}g_u(t)=\sum_{i=1}^n\frac{1}{\sigma^2_{u,i}(t)}=\sum_{i=1}^n(1-h_i(t_0+t)/u)^2.
 \EQN
 Then by \eqref{eq:gtt0} and the fact that $\min_{1\leq i\leq n}c_i>0$, we have  for $\theta>0 $ sufficiently small and $u$ sufficiently large
\BQNY
g_u(t)-g_u(0)&=&\sum_{i=1}^n(1-h_i(t_0+t)/u)^2-\sum_{i=1}^n(1-h_i(t_0)/u)^2\\
&\geq & \sum_{i=1}^n\frac{h_i(t_0)-h_i(t_0+t)}{u}\\
&\geq & \mathbb{C}_1\frac{|t|^\gamma}{u} \geq \mathbb{C}_1\frac{(\ln u)^q}{u^2}, \quad t\in [t_0-\theta, t_0+\theta] \setminus  (t_0+\E).
 \EQNY
 Consequently, there exists $C>0$ such that
\BQNY
\sup_{t\in [t_0-\theta, t_0+\theta] \setminus  (t_0+\E)}\Var(Y_u(t))
=\sup_{t\in [t_0-\theta, t_0+\theta] \setminus  (t_0+\E)}\frac{1}{g_u(t)}\leq\frac{1}{g_u(0)+\frac{C(\ln u)^q}{u^2}}.
\EQNY
Moreover, for $\theta>0 $ sufficiently small and $u$ sufficiently large
\BQN\label{gt1}
g_u(t)-g_u(0)
&\geq &\frac{ \sum_{i=1}^n h_i(t_0)- \sum_{i=1}^n h_i(t_0+t)}{u}\nonumber\\
&\geq & \frac{\mathbb{C}_2}{u}, \quad t\in [0,T]\setminus[t_0-\theta, t_0+\theta].
 \EQN
 Thus there exists $C_1>0$ such that \BQN\label{upper}
\sup_{t\in [0,T]\setminus[t_0-\theta, t_0+\theta]}\Var(Y_u(t))
=\sup_{t\in [0,T]\setminus[t_0-\theta, t_0+\theta]}\frac{1}{g_u(t)}\leq\frac{1}{g_u(0)+\frac{C_1}{u}}.
\EQN
Consequently,
\BQNY
\sup_{t\in [t_0-\theta, t_0+\theta] \setminus  (t_0+\E)}\Var(Y_u(t))
\leq\frac{1}{g_u(0)+\frac{C_2(\ln u)^q}{u^2}}, \quad [0,T]\setminus (t_0+\E),
\EQNY
with $C_2>0$.
Moreover,
in light of (\ref{holderf}) and (\ref{stationaryR0}), we have that
\BQN\label{holderyy}
\mathbb{E}\LT(Y_u(t)-Y_u(s)\RT)^2
\leq\mathbb{C}_3|t-s|^{\mu}, \quad s,t\in [0,T]
\EQN
for $\mu>0$. Piterbarg inequality (Theorem 8.1 in \cite{Pit96}) leads to
\BQNY
\Pi_1(u)&\leq&\pk{\sup_{t\in [0,T] \setminus  (t_0+\E)} Y_u(t)>u}\\
&\leq& \mathbb{C}_3 u^{2/\mu} \Psi\left(u\sqrt{g_u(0)+\frac{C_2(\ln u)^q}{u^2}}\right)\\
&=& o(\Pi(u)), \quad u\rw\IF.
\EQNY
This establishes the claim.\\
ii)  Without loss of generality, we assume that $0<A<B<T$. Then for $\epsilon>0$ sufficiently small
\BQNY
\pk{\exists_{t\in[A,B]} \LT(\vk{X}(t)+{\vk{h}}(A)\RT)  >u\vk{1}}&\le& \pk{\exists_{t\in[0,T]} \LT(\vk{X}(t)+{\vk{h}}(t)\RT)  >u\vk{1}}\\
&\le&\pk{\exists_{t\in[0,A-\epsilon]} \LT(\vk{X}(t)+{\vk{h}}(t)\RT)  >u\vk{1}}+\pk{\exists_{t\in[A-\epsilon, B+\epsilon]} \LT(\vk{X}(t)+{\vk{h}}(A)\RT) >u\vk{1}}\\
&&+\pk{\exists_{t\in[B+\epsilon,T]} \LT(\vk{X}(t)+{\vk{h}}(t)\RT)  >u\vk{1}}.
\EQNY

In view of (\ref{stationaryR0}) and (\ref{stationaryR2}) and by Theorem 4.1 in  \cite{Tabis}, we have that for any $0\leq x<y\leq T$
\BQNY
\pk{\exists_{t\in[x,y]} \LT(\vk{X}(t)+{\vk{h}}(A)\RT) >u\vk{1}}&=&\pk{\exists_{t\in[x,y]}\vk{X}(t)>u\vk{1}-\vk{h}(A)}\\
&\sim &u^{\frac{2}{\alpha}}\int_{x}^{y}\mathcal{H}_{\alpha,\vk{a}(t)
\mathbf{I}_{\{\vk{\alpha}=\alpha\vk{1}\}}}dt\prod_{i=1}^n\Psi\LT(u-h_{m,i}\RT),\quad u\rw\IF,
\EQNY
where  $\int_{x}^{y}\mathcal{H}_{\alpha,\vk{a}(t)\mathbf{I}_{\{\vk{\alpha}=\alpha\vk{1}\}}}dt$ is a finite and positive constant (see \cite{Tabis}).
Next we show that $\pk{\exists_{t\in[0,A-\epsilon]} \LT(\vk{X}(t)+{\vk{h}}(t)\RT)  >u\vk{1}}$ is  negligible. Rewrite
\BQNY
\pk{\exists_{t\in[0,A-\epsilon]} \LT(\vk{X}(t)+{\vk{h}}(t)\RT)  >u\vk{1}}=\pk{\exists_{t\in[0,A-\epsilon]} Y_u(t)>u},
\EQNY
where $Y_u$ is defined in (\ref{Yu}). Note  that (\ref{upper}) still holds in the case considered with $[0,A-\epsilon]$ instead of  $[0,T]\setminus[t_0-\theta, t_0+\theta]$. Therefore, in view of (\ref{holderyy}), by Piterbarg inequality we have that
$$\pk{\exists_{t\in[0,A-\epsilon]} Y_u(t)>u}\leq  \mathbb{C}_4 u^{2/\mu} \Psi\left(u\sqrt{g_u(0)+\frac{C_1}{u}}\right)=o\left(u^{\frac{2}{\alpha}}
\prod_{i=1}^n\Psi\LT(u-h_{m,i}\RT)\right), \quad u\rw\IF.$$
Analogously,
$$\pk{\exists_{t\in[B+\epsilon,T]} \LT(\vk{X}(t)+{\vk{h}}(t)\RT)  >u\vk{1}}=o\left(u^{\frac{2}{\alpha}}\prod_{i=1}^n\Psi\LT(u-h_{m,i}\RT)\right), \quad u\rw\IF.$$
Therefore, we conclude that as $u\rw\IF$
\BQNY
u^{\frac{2}{\alpha}}\int_{A}^{B}\mathcal{H}_{\alpha,\vk{a}(t)
\mathbf{I}_{\{\vk{\alpha}=\alpha\vk{1}\}}}dt\prod_{i=1}^n\Psi\LT(u-h_{m,i}\RT) &\leq& \pk{\exists_{t\in[0,T]} \LT(\vk{X}(t)+{\vk{h}}(t)\RT)  >u\vk{1}}\\
&\leq& u^{\frac{2}{\alpha}}\int_{A-\epsilon}^{B+\epsilon}\mathcal{H}_{\alpha,\vk{a}(t)
\mathbf{I}_{\{\vk{\alpha}=\alpha\vk{1}\}}}dt
\prod_{i=1}^n\Psi\LT(u-h_{m,i}\RT).
\EQNY
We establish the claim by letting $\epsilon\rw 0$ in the above inequalities. This completes the proof.
\QED

\proofprop{prop1}
We notice that
\begin{align*}
p(u)=\pk{\exists_{t\in[0,T]}\LT(\vk{B}_{\vk{\alpha}}(t)-\vk{c}t\RT) >u\vk{d}}
=\pk{\exists_{t\in[0,T]}\LT(\frac{1}{\vk{d}}\vk{B}_{\vk{\alpha}}(t)-\frac{\vk{c}t}{\vk{d}}\RT) >u\vk{1}},
\end{align*}
and the variance function $\sigma^2_i(t)$ and correlation function $r_i(s,t)$ of $\frac{B_{\alpha_i}(t)}{d_i}$  satisfy
\begin{align*}
&r_i(s,t)=1-\frac{1}{2T^{\alpha_i}}|t-s|^{\alpha_i}(1+o(1)), s,t\rw T,\\
& \sigma_i(t)=\frac{T^{{\alpha_i}/2}}{d_io}-\frac{\alpha_i}{2d_i}T^{{\alpha_i}/2-1}(T-t)(1+o(1)), t\rw T,
\end{align*}
where $T$ is the unique maximum point of $\sigma_i(t), 1\leq i\leq n$ over $[0,T]$. Moreover,
$$-\frac{c_it}{d_i}=-\frac{c_iT}{d_i}+\frac{c_i}{d_i}|T-t|, \quad t\rw T.$$
Therefore, in light of \netheo{Thm2} and \nekorr{Korr1}, we have that
\begin{align*}
\pk{\exists_{t\in[0,T]}\LT(\vk{B}_{\vk{\alpha}}(t)-\vk{c}t\RT) >u\vk{d}}\sim u^{(\frac{2}{\alpha}-2)_{+}}\prod_{i=1}^n\Psi\LT(\frac{d_i u+c_iT}{T^{{\alpha_i}/2}}\RT)
\LT\{
\begin{array}{ll}
\mathcal{H}_{\alpha,\vk{\varsigma}\mathbf{I}_{\{\vk{\alpha}=\alpha\vk{1}\}}}
\int_{0}^\IF e^{-\sum_{i=1}^n f_i(t)}dt,&\ \ \text{if}\ \alpha<1,\\
\mathcal{P}^{\vk{f}}_{\alpha,\vk{\varsigma}
\mathbf{I}_{\{\vk{\alpha}=\alpha\vk{1}\}}}[0,\IF), &\ \ \text{if}\ \alpha=1,\\
1,&\ \ \text{if}\ \alpha>1,
\end{array}
\RT.
\end{align*}
and
\BQNY
\pk{(T-\tau_u)u^2\leq x\Big| \tau_u\leq T}\sim
\LT\{
\begin{array}{ll}
1-e^{-\LT(\sum_{i=1}^n\frac{\alpha_id_i^2}{2T^{\alpha_i+1}}\RT)x},&\ \ \text{if}\ \alpha<1,\\
{\mathcal{P}^{\vk{f}}_{\alpha,\vk{\varsigma}
\mathbf{I}_{\{\vk{\alpha}=\alpha\vk{1}\}}}[0,x]}\Big/
{\mathcal{P}^{\vk{f}}_{\alpha,\vk{\varsigma}
\mathbf{I}_{\{\vk{\alpha}=\alpha\vk{1}\}}}[0,\IF)}, &\ \ \text{if}\ \alpha=1,\\
1,&\ \ \text{if}\ \alpha>1,
\end{array}
\RT.
\EQNY
where $\alpha=\min_{1\leq i\leq n}\alpha_i$, \bl{$\vk{\varsigma}=(\varsigma_1,\dots, \varsigma_n)$ with $\varsigma_i=\frac{d_i^2}{2T^{2\alpha_i}}$ }and $f_i(t)=\frac{\alpha_id_i^2}{2T^{\alpha_i+1}}|t|$.
\QED

\section{Appendix}	
\proofkorr{Korr1}
By definition,
\BQN\label{korr1eq1}
\pk{(T-\tau_u)u^{2/\beta}\leq x\Big| \tau_u\leq T}
=\frac{\pk{\exists_{t\in[T-u^{-2/\beta}x,T]}\LT(\vk{X}(t)+{\vk{h}}(t)\RT) >u\vk{1}}}{\pk{\exists_{t\in[0,T]}\LT(\vk{X}(t)+{\vk{h}}(t)\RT) >u\vk{1}}}
\EQN
The asymptotics of denominator in (\ref{korr1eq1}) follows
by  \netheo{Thm2}.
In order to get the asymptotics of nominator of (\ref{korr1eq1})
we follow the same argument as in the proof of \netheo{Thm2} (part related with the asymptotics of $\Pi_1(u)$),
which leads to
% with $E(u)=[-u^{-2/\beta}x, 0]$, we have with the same natation there
\BQN\label{korr1eq2}
\pk{\exists_{t\in[T-u^{-2/\beta}x, T]}\LT(\vk{X}(t)+{\vk{h}}(t)\RT) >u\vk{1}}&\sim& u^{(\frac{2}{\alpha}-\frac{2}{\beta})_{+}}\prod_{i=1}^n\Psi\LT(\frac{u-h_i(t_0)}
{\sigma_i(t_0)}\RT)
\nonumber\\
&&\times\LT\{
\begin{array}{ll}
\mathcal{H}_{\alpha,\frac{\vk{a}}{\vk{\sigma}^2(t_0)}
\mathbf{I}_{\{\vk{\alpha}=\alpha\vk{1}\}}}\int_{-x}^0 e^{-\sum_{i=1}^n f_i(x)}dx,&\ \ \text{if}\ \alpha<\beta,\\
\mathcal{P}^{\vk{f}}_{\alpha,\frac{\vk{a}}{\vk{\sigma}^2(t_0)}
\mathbf{I}_{\{\vk{\alpha}=\alpha\vk{1}\}}}[-x,0], &\ \ \text{if}\ \alpha=\beta,\\
1,&\ \ \text{if}\ \alpha>\beta,
\end{array}
\RT.
\EQN
which completes the proof.
%Inserting \eqref{korr1eq2} and the result of \netheo{Thm2} about $\pk{\exists_{t\in[0,T]}\LT(\vk{X}(t)+{\vk{h}}(t)\RT) >u\vk{1}}$ into \eqref{korr1eq1}, we obtain our claim.
\QED

\COM{{\bf{Proof of (\ref{e.2}):}}
%\proofkorr{ex1ruin}
By definition,
\BQN\label{ratio}
\pk{(T-\tau_u)u^2\leq x\Big| \tau_u\leq T}=\frac{\pk{\exists_{t\in[T-xu^{-2},T]} \LT(\vk{B}_{\vk{\alpha}}(t)+\vk{c}t\RT) >u\vk{d}}}{\pk{\exists_{t\in[0,T]} \LT(\vk{B}_{\vk{\alpha}}(t)+\vk{c}t\RT) >u\vk{d}}}.
\EQN
For all $u>0$ large,
\BQNY
\pk{\exists_{t\in[T-xu^{-2},T]} \LT(\vk{B}_{\vk{\alpha}}(t)+\vk{c}t\RT) >u\vk{d}}
&=&\pk{\exists_{t\in[-xu^{-2},0]} \LT(\vk{B}_{\vk{\alpha}}(T+t)+\vk{c}(T+t)\RT) >u\vk{d}}\\
&=&\pk{\exists_{t\in[-xu^{-2},0]}\frac{ \vk{B}_{\vk{\alpha}}(T+t)}{\vk{d}-\vk{c}(T+t)/(u\vk{1})} >u\vk{1}}
\EQNY
Denote by $\vk{X}_u(t)=\frac{ \vk{B}_{\vk{\alpha}}(T+t)}{\vk{d}-\vk{c}(T+t)/(u\vk{1})}$.
Let $\sigma_{u,i}(t)$ and $r_{u,i}(s,t)$ denote the variance and correlation functions of $X_{u,i}(t)$, respectively. Next we check the conditions of \netheo{PreThm1}. {\bf A1} holds straightforwardly. Direct calculation gives that
 $$
 \lim_{u\rightarrow\IF}{\sup_{t\in [-xu^{-2},0]}}\LT|\frac{\LT(\frac{\sigma_{u,i}(0)}{\sigma_{u,i}(t)}-1\RT)u^2-f
 (u^2t)}
 {\LT|f(u^2t)\RT|+1}\RT|=0,
$$
 where $f(t)=\frac{\alpha}{2T}|t|$, and
$$
\lim_{u\rightarrow\IF}\sup_{s,t\in [-xu^{-2},0]}\left|\frac{1-r_{u,i}(t,s)}{2^{-1}T^{-\alpha_i}|t-s|^{\alpha_i}}-1\right|=0.
$$
These confirm {\bf A2-A3} hold.
Moreover, it easy to check that (\ref{xx}) is satisfied with $x_1=-x$ and $x_2=0$. Consequently, in light of  \netheo{PreThm1}, we have that as $u\rw\IF$
\BQNY
\pk{\exists_{t\in[T-xu^{-2},T]} \LT(\vk{B}_{\vk{\alpha}}(t)+\vk{c}t\RT) >u\vk{d}}\sim u^{(\frac{2}{\alpha}-2)_{+}}\prod_{i=1}^n\Psi\LT(\frac{d_i u-c_iT}{T^{{\alpha_i}/2}}\RT)
\LT\{
\begin{array}{ll}
\mathcal{H}_{\alpha, \vk{\varsigma}\mathbf{I}_{\{\vk{\alpha}=\alpha\vk{1}\}}}
\int_{-x}^0 e^{-\sum_{i=1}^n f_i(t)}dt,&\ \ \text{if}\ \alpha<1,\\
\mathcal{P}^{\vk{f}}_{\alpha, \vk{\varsigma}\mathbf{I}_{\{\vk{\alpha}
=\alpha\vk{1}\}}}[-x,0], &\ \ \text{if}\ \alpha=1,\\
1,&\ \ \text{if}\ \alpha>1,
\end{array}
\RT.
\EQNY
with the same $\alpha, \vk{\varsigma}$ and $\vk{f}$ as in the Proof of (\ref{e.1}).
Inserting the above asymptotics into (\ref{ratio}), we establish the claim.
\QED}

\section*{Acknowledgments}
Thanks to  the Swiss National Science Foundation Grant 200021-175752/1, whereas K. D\c ebicki
also acknowledges partial support from NCN Grant No 2015/17/B/ST1/01102 (2016-2019).

 	\bibliographystyle{ieeetr}

\bibliography{VecValue}
\end{document}